\documentclass{statsoc}
 
\usepackage[a4paper]{geometry}
\usepackage{graphicx}
\usepackage[textwidth=8em,textsize=small]{todonotes}
\usepackage{amsmath} 
\usepackage{natbib}

\RequirePackage[OT1]{fontenc}
\RequirePackage{amssymb,lmodern,epstopdf,enumitem,lmodern} 
\RequirePackage[colorlinks,citecolor=blue,urlcolor=blue]{hyperref}
 
\usepackage{bm,color,fancyvrb,xcolor,mathtools}
\usepackage{amscd} 
\usepackage{mathrsfs}
\usepackage{bbm}
\usepackage{array}
\usepackage{multirow,verbatim} 
\usepackage{diagbox,array}

\def\eps{\varepsilon}
\font\tencmmib=cmmib10 \skewchar\tencmmib '60
\newfam\cmmibfam
\textfont\cmmibfam=\tencmmib

\def\lessim{\ \lower4pt\hbox{$
		\buildrel{\displaystyle <}\over\sim$}\ }
\def\gessim{\ \lower4pt\hbox{$\buildrel{\displaystyle >}
		\over\sim$}\ }

\usepackage{fancyhdr}
\usepackage{algorithm,algorithmicx}%,algpseudocode}
\usepackage[noend]{algpseudocode}
\usepackage{tensor}
\usepackage{verbatim}

\newtheorem{theorem}{Theorem}[section]
\newtheorem{proposition}[theorem]{Proposition}
\newtheorem{lemma}{Lemma}
\newtheorem{defn}[theorem]{Definition}
\newtheorem{corollary}[theorem]{Corollary}
\newtheorem{remark}{Remark}[section]

\DeclarePairedDelimiter\floor{\lfloor}{\rfloor}

\DeclareFontFamily{U}{mathx}{\hyphenchar\font45}
\DeclareFontShape{U}{mathx}{m}{n}{<-> mathx10}{}
\DeclareSymbolFont{mathx}{U}{mathx}{m}{n}
\DeclareMathAccent{\widebar}{0}{mathx}{"73}

\renewcommand{\hat}{\widehat}
\renewcommand{\tilde}{\widetilde}

\renewcommand{\hat}{\widehat}

\newcommand{\bbm}[1]{\text{\usefont{U}{bbm}{m}{n}#1}}
\newcommand\numberthis{\addtocounter{equation}{1}\tag{\theequation}}

     \def\EE{\mathbb{E}}

\def\bm{\bfm m}

     \def\PP{\mathbb{P}}
     
     \def\RR{\mathbb{R}}

 \def\cA{{\cal  A}}

\def\calE{{\cal  E}} \def\cE{{\cal  E}}

\def\calH{{\cal  H}}

\def\calK{{\cal  K}} 
 
\def\calM{{\cal  M}} \def\cM{{\cal  M}}

\def\calS{{\cal  S}} \def\cS{{\cal  S}}
 \def\cT{{\cal  T}}

\def\calW{{\cal  W}}

\def\calZ{{\cal  Z}} \def\cZ{{\cal  Z}}

\def\eps{\varepsilon}

\newcommand{\vertiii}[1]{{\left\vert\kern-0.25ex\left\vert\kern-0.25ex\left\vert #1 
	\right\vert\kern-0.25ex\right\vert\kern-0.25ex\right\vert}}

\usepackage{tikz}
\usetikzlibrary{decorations.pathreplacing,calc}

\def\scrT{\mathscr{T}} \def\scrZ{\mathscr{Z}}

\renewcommand{\bm}{\boldsymbol}
\def\P{\mathbb{P}}

\title[Minimax Optimal Transfer Learning for Classification under Distributed Differential Privacy]{Minimax And Adaptive Transfer Learning for Nonparametric Classification under Distributed Differential Privacy Constraints}
\author[Auddy, Cai, and Chakraborty]{Arnab Auddy}
\address{Department of Biostatistics, Epidemiology and Informatics,\\ University of Pennsylvania, Philadelphia, PA 19104.}
\author[Auddy, Cai, and Chakraborty]{T. Tony Cai}
\address{Department of Statistics and Data Science, The Wharton School,\\ University of Pennsylvania, Philadelphia, PA 19104.}
\author[Auddy, Cai, and Chakraborty]{Abhinav Chakraborty}
\address{Department of Statistics and Data Science, The Wharton School,\\ University of Pennsylvania, Philadelphia, PA 19104.}

\begin{document}
	\begin{abstract}
		This paper considers minimax and adaptive transfer learning for nonparametric classification under the posterior drift model with distributed differential privacy constraints. Our study is conducted within a heterogeneous framework, encompassing diverse sample sizes, varying privacy parameters, and data heterogeneity across different servers. 
		
		We first establish the minimax misclassification rate, precisely characterizing the effects of privacy constraints, source samples, and target samples on classification accuracy. The results reveal interesting phase transition phenomena and highlight the intricate trade-offs between preserving privacy and achieving classification accuracy. We then develop a data-driven adaptive classifier that achieves the optimal rate within a logarithmic factor across a large collection of parameter spaces while satisfying the same set of differential privacy constraints. Simulation studies and real-world data applications further elucidate the theoretical analysis with numerical results.
	\end{abstract}

	%%%%%%%%%%%%%%%%%%%%%%%%%%%%%%%%%%%%%%%%%%
	\section{Introduction}
	%%%%%%%%%%%%%%%%%%%%%%%%%%%%%%%%%%%%%%%%%%

	Massive and diverse datasets are now routinely collected across a wide range of scientific fields, including genomics, neuroimaging, astrophysics, climate studies, and signal processing. In many applications, alongside the primary data from the target study, additional datasets from different populations or environments with similar structures have also been collected. Transfer learning, which aims to improve learning performance in a target domain by transferring knowledge from different but related source domains, has become a vibrant and promising area of research in machine learning. This concept has found applications in areas such as computer vision, speech recognition, and genre classification.
	
%	In many scenarios, labeled data is available from a distribution $P$ (the source distribution), while a relatively small quantity of labeled or unlabeled data is drawn from a distribution $Q$ (the target distribution). The goal is to make statistical inference under $Q$, leveraging information from $P$. Several approaches have been proposed to measure the similarity between $P$ and $Q$, including divergence bounds, covariate shift, label shift, and posterior drift.
	 
	Alongside the transfer learning framework, another crucial consideration in modern data science is the preservation of privacy. Sensitive data are often spread across various sources, each presenting unique challenges in privacy preservation (see, e.g., \cite{guo2024comprehensive}). Differential privacy (DP) has emerged as a leading framework for ensuring that statistical analysis results do not compromise the confidentiality of individual data points. Originally introduced by \cite{dwork2006calibrating}, DP has garnered significant academic attention and has been embraced by industry leaders like Google, Microsoft, and Apple, as well as governmental entities such as the US Census Bureau (\cite{abowd2016challenge}).
	
	Addressing the distributed nature of data collection and analysis is crucial due to its implications for privacy preservation and collaboration. This raises a number of natural questions: firstly, when and how can inference under the target distribution be improved by leveraging data points distributed across several sources with disparate privacy constraints? Secondly, what are the fundamental limits of inference for the target with privacy-constrained learning from the sources?
	
	The two questions outlined above frequently arise in modern data analysis, spurring a flurry of recent research. The answers to both questions fundamentally depend on the specific inference task at hand. Only recently has rigorous research explored the theoretical performance of transfer learning  under differential privacy constraints, addressing both various parametric problems (\cite{li2024federated}) and nonparametric problems (\cite{ma2023optimal}).

	In this paper, we focus on the task of binary nonparametric classification, also referred to as domain adaptation in the literature. Our objective is to achieve statistically optimal transfer learning under distributed DP constraints. We tackle the challenge of heterogeneous data sources with distinct distributions, examining how privacy parameters, sample sizes, and data heterogeneity affect classification performance. Our proposed classifiers are designed to adapt to unknown data heterogeneity and parameters, while maintaining statistical optimality. This adaptability enhances performance even under strict privacy constraints, balancing privacy preservation with classification accuracy. Before delving into further details, we outline some fundamental concepts of classification with transfer learning.
	
	In a classification problem from a single source, we observe independent copies of a tuple $(X,Y)$ from a distribution $P$ where $X\in\RR^d$ are the covariates, and $Y\in\{0,1\}$ denotes the binary class labels. To measure the efficacy of transferring information from a source distribution $P$ to a different target distribution $Q$, several methods have been proposed to quantify the similarity of $P$ and $Q$. Building upon the intuition that transfer is easier if both $P$ and $Q$ have high masses in a common region, divergence measures have been used on either the covariate space or as discrepancy between labels. See, e.g., \cite{ben2010theory, germain2013pac, cortes2019adaptation, sugiyama2007direct} and references therein. While these are general metrics of measuring the difference between $P$ and $Q$, such approaches often tend to be pessimistic, as shown by \cite{kpotufe2021marginal}. Incorporating the structure of the classification problem leads to describing more specific transfer models such as covariate shift: where one posits a difference in the marginal distribution of $X$ from $P$ to $Q$, but assumes the posterior probabilities of $\PP(Y=1|X)$ to be the same for both $P$ and $Q$. On the other hand, label shift assumes that the class probability $\PP(Y=1)$ differs from $P$ to $Q$, but the class-conditional covariate distributions $\PP(X|Y=1)$ remains the same. For more details we refer the reader to \cite{kpotufe2021marginal,sugiyama2007direct} for the covariate shift, and to \cite{garg2020unified,lipton2018detecting,maity2022minimax} for the label shift paradigms. More flexible transfer mechanisms have also been considered, see, e.g., \cite{reeve2021adaptive, fan2023robust}.
 	
	In this paper we focus on the posterior drift model, where the covariates have the same marginal distribution under both the source and target distributions, but the posterior probabilities $\PP(Y=1|X)$ undergoes a shift from the source to the target. Posterior drift has often been studied in the literature: see \cite{cai2021transfer,liu2020computationally,maity2024linear,scott2019generalized} and references therein. The posterior drift model naturally arises in a range of applications where the data is distributed and there are privacy concerns. Here are a few examples:
	\begin{itemize}
		\item \textbf{Healthcare Monitoring Across Hospitals:} In applications  where multiple hospitals contribute data for healthcare monitoring (\cite{yeung2019local}) (e.g., patient vital signs, medical histories), each hospital's data is sensitive and subject to privacy regulations like HIPAA. As medical practices evolve and patient demographics change, the underlying distribution of health data in each hospital's domain may drift over time. However, due to privacy concerns, the data cannot be aggregated into a central repository for analysis (see for e.g., \cite{ju2020federated} which proposes a privacy-preserving federated transfer learning architecture for EEG classification, achieving higher accuracy without data sharing.). Consequently, the transfer learning model trained on data from one hospital may experience posterior drift when applied to another hospital's data, leading to degradation in performance over time.
		\item \textbf{Financial Fraud Detection in Banking Networks:} Banks collaborate to detect financial fraud by sharing transaction data (\cite{phua2010comprehensive}, \cite{809570}), but due to privacy regulations and competitive concerns, they cannot directly share sensitive customer information. Over time, patterns of financial fraud may evolve due to changes in customer behavior, economic conditions, or fraud tactics. However, because of customer privacy concerns, each bank must maintain control over its own data, making it challenging to aggregate data for analysis. As a result, transfer learning models used for fraud detection (\cite{lebichot2020deep}) may experience posterior drift as the data distributions in different banks' domains shift over time.  
		\item \textbf{Social Media Analysis Across Platforms:} Social media platforms collect user-generated content and engagement metrics for analyzing trends, sentiment analysis, and targeted advertising (\cite{saura2019comparing}). However, due to privacy regulations and platform policies, individual user data cannot be shared openly between social media platforms. As online communities evolve and user behaviors shift with trending topics, viral content, and platform updates, the underlying distribution of social media data in each platform's domain may vary. Therefore, learning methods for social media analysis (\cite{wang2020cross})  may encounter posterior drift when applied across platforms, while simultaneously requiring to protect sensitive user data. 
	\end{itemize}  
	
	Nonparametric classification in the posterior drift model has been considered previously by \cite{cai2021transfer}, where the minimax rate of misclassification risk without the privacy constraints is established. The present work builds upon their framework, and establishes the optimal rates for classification using data that are distributed across servers with varying quality and different privacy constraints. Our results reveal interesting phase transition phenomena and highlight the intricate trade-offs between preserving privacy and achieving classification accuracy. We propose statistically optimal adaptive procedures that effectively balances this trade-off between privacy and accuracy.
	
	Suppose the source and target distributions $P$ and $Q$ have similar covariate distributions, and consider the functions $\eta_P(x):=\PP(Y=1|X=x)$ under $P$, and $\eta_Q(x):=\PP(Y=1|X=x)$ under $Q$. Further, suppose that both $P$ and $Q$ have the same decision boundary, that is, $(\eta_Q(x)-1/2)(\eta_P(x)-1/2)\ge 0$ for all $x\in\RR^d$. In this framework, \cite{cai2021transfer} quantified the efficiency of transfer learning at $Q$ using data from $P$ through the so-called relative signal exponent $\gamma$, where $\gamma>0$ is a number such that 
	$$
	\left\vert\eta_P(x)-\frac{1}{2}\right\vert 
	\ge \left\vert\eta_Q(x)-\frac{1}{2}\right\vert^{\gamma}.
	$$
	When $\gamma<1$, the covariates are better separated into the two classes in the source distribution $P$ than in $Q$. In such a case, the $P$-data are especially useful for transfer learning. The situation reverses when $\gamma>1$. In this paper, we describe a kernel based classifier that leverages the information from sources, while preserving data privacy in the DP framework. For both the source and the target, we compute:
	\begin{equation}\label{eq:kernel-estimators}
	T_P(x)=\dfrac{1}{n_Ph^d}\sum_{i=1}^{n_P}\left(Y_i-\frac{1}{2}\right)K\left(\frac{X_i-x}{h}\right)
	\text{ and }
	T_Q(x)=\dfrac{1}{n_Qh^d}\sum_{i=1}^{n_Q}\left(Y_i-\frac{1}{2}\right)K\left(\frac{X_i-x}{h}\right),
	\end{equation}
	for a suitable kernel $K(\cdot)$ and bandwidth $h$. It can be checked that $T_P(x)$ and $T_Q(x)$ are pointwise consistent estimators of $\left(\eta_P(x)-\frac{1}{2}\right)g_P(x)$ and $\left(\eta_Q(x)-\frac{1}{2}\right)g_Q(x)$ respectively, where $g_P(\cdot)$ and $g_Q(\cdot)$ are the joint densities of the covariates under $P$ and $Q$ respectively. Thus the sign of a suitable linear combination of $T_P(x)$ and $T_Q(x)$ defines a reasonable transfer learning classifier. 
	
	The above classification procedure is however, not differentially private. Nonparametric classifiers with local privacy have previously been considered by \cite{berrett2019classification, ma2023optimal}. In this paper, we focus on ``server-level" privacy, where each server can access its own set of unperturbed data, but imposes privacy constraints when sharing aggregated information with other servers. Due to its immense applicability in practice, such a distributed privacy setting has recently received significant attention in the literature: see e.g., \cite{Cai2023FL-NP-Regression,acharya2023discrete,liu2020learning,levy2021learning}. 
	
	To emphasize the impact of distributed data privacy, we consider a scenario with $m$ source servers from distribution $P$, each holding $n_P$ samples. Similar to~\eqref{eq:kernel-estimators} we define the estimator for each $j$-th source server, denoted as $T^{(j)}_P(\cdot)$. For attaining differential privacy, it is standard to add external noise, with an important tradeoff in mind. The noise level must be enough to ensure the $(\eps,\delta)$-differential privacy guarantees (see Definition~\ref{def:diff-privacy}) for both the source and the target, but it has to be added prudently so as to not worsen the classification accuracy too much. In our case, we output a noisy test function as follows. We use the Gaussian process based mechanism used by \cite{hall2013differential} to privatize kernel estimators. With $T_Q(\cdot)$ and $T^{(j)}_P(\cdot)$ as defined above, we compute:
	$$
	\tilde{T}_Q(x)=T_Q(x)+\sigma_Q\xi_Q(x)
	\text{ and }
	\tilde{T}^{(j)}_P(x)=T^{(j)}_P(x)+\sigma_P \xi^{(j)}_P(x)\,\text{ for }\,j=1,\ldots,m,
	$$
	where $\xi_P(\cdot)$ and $\xi^{(j)}_Q(\cdot)$ are Gaussian processes with covariance function $K(\cdot/h)$, while $\sigma_P^2$ and $\sigma_Q^2$ are suitably chosen noise variances that depend on the privacy parameters $(\eps,\delta)$, sample sizes and the bandwidth $h$. Next, we choose an appropriate weight $w\in[0,1]$, and define the privatized transfer weighted estimator 
	\[
	\tilde{T}(x)=w\tilde{T}_Q(x)+(1-w)\left(\frac{1}{m}\sum_{j=1}^m\tilde{T}^{(j)}_P(x)\right).
	\]
	Finally our classifier becomes:
	\begin{equation}\label{eq:intro-classifier}
		\hat{f}(x):=\bbm{1}(\tilde{T}(x)\ge 0).
	\end{equation}
	We use the standard notion of excess risk (see, e.g., \cite{audibert2007plugin}) to evaluate the performance of our classifier. Let $f_Q^*(x)=\bbm{1}(\eta_Q(x)\ge 1/2)$ be the Bayes classifier for the target distribution. Then the excess risk of a classifier under the target distribution is defined as:
	$$
	\calE_Q(\hat{f})=\EE\left[\PP_Q(Y\neq \hat{f})\right] - \PP_Q(Y\neq f^*_Q).
	$$
	We assume that the margin assumption (see Definition~\ref{def:margin-assn}) holds with parameter $\alpha$ for $Q$, and $\eta_Q(x)$ is $(\beta,L_{\beta})$-H{\"o}lder for some $0<\beta<1$. Then if the target and the source are both constrained to be $(\eps,\delta)$ differentially private, we prove the following minimax rate:
	\begin{align*}\label{eq:minimax-rate-intro}
	&~\inf_{\tilde{f}\in\calM(\eps,\delta)}
	\sup_{P,Q}
	\calE_Q(\hat{f})\\
	\asymp&~
	\left[ L_N\left\{\left( n_Q^{-\frac{\beta(1+\alpha)}{2\beta + d}} \vee (n_Q^2\varepsilon^2)^{-\frac{\beta(1+\alpha)}{2\beta + 2d}}\right) \bigwedge \left( (mn_P)^{-\frac{\beta(1+\alpha)}{2\beta \gamma+ d}} \vee (mn_P^2\varepsilon^2)^{-\frac{\beta(1+\alpha)}{2\beta \gamma+ 2d}}\right) \right\}\wedge 1\right]\numberthis
	\end{align*} 
	where $L_N$ is of order at most ${\rm polylog}(mn_P+n_Q)$, the infimum is over all possible transfer learning classifiers satisfying $(\eps,\delta)$-DP, and the supremum is over all distributions in the posterior drift framework. The minimax rate is attained by the classifier $\hat{f}$ in \eqref{eq:intro-classifier}. 
%	The above rate is in fact nearly optimal, since we also prove the following lower bound:
%	\begin{align*}
%		&~\inf_{\tilde{f}\in\calM(\eps,\delta)}
%		\sup_{P,Q}
%		\calE_Q(\tilde{f})\\
%		\gtrsim&~
%		\left[\left\{\left( n_Q^{-\frac{\beta(1+\alpha)}{2\beta + d}} \vee (n_Q^2\varepsilon^2)^{-\frac{\beta(1+\alpha)}{2\beta + 2d}}\right) \bigwedge \left( (mn_P)^{-\frac{\beta(1+\alpha)}{2\beta \gamma+ d}} \vee (mn_P^2\varepsilon^2)^{-\frac{\beta(1+\alpha)}{2\beta \gamma+ 2d}}\right) \right\}\wedge 1\right]
%	\end{align*}
	
\begin{figure}[h]
		\centering
		\includegraphics[scale = 0.65]{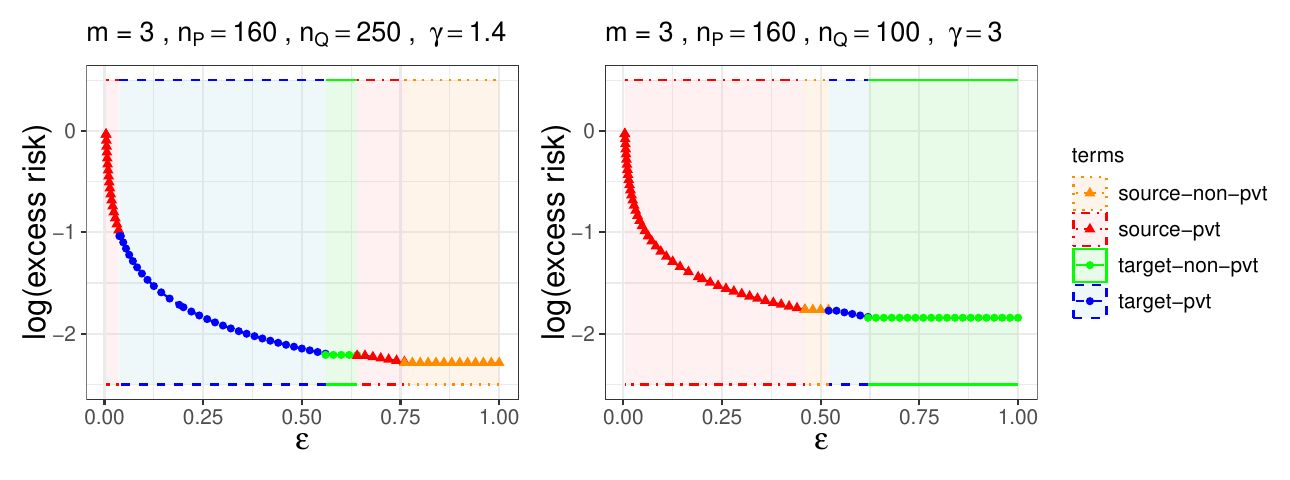}
		\caption{Relationship of logarithm of excess risk with $\varepsilon$ as given by~\eqref{eq:minimax-rate-homogenous}, the smoothness level $\beta = 0.25$ and dimension $d =2$.
		 }
		\label{fig:theo-rates}
	\end{figure}

	The minimax rate is determined by a trade-off between four quantities: the non-private rates for the source and the target, as well as their privatized counterparts. An interesting phenomenon emerges through the relative effect of the efficacy of transfer and the cost of privacy, characterized by the parameters $\gamma$ and $\eps$ along with server size $m$ and sample sizes $n_P, n_Q$. This discrepancy essentially means that a higher privacy cost is incurred when data is more distributed, see Figure~\ref{fig:theo-rates}. For a detailed discussion, refer to the remarks following Theorem \ref{th:homogeneous-minimax-rate}. 
	
	Figure~\ref{fig:theo-rates} illustrates the minimax rate as described in Equation \eqref{eq:minimax-rate-intro}, showcasing the different regimes of differentially private transfer learning excess risk. Several intriguing phenomena are observed. For instance, across varying privacy budgets $\varepsilon$, the regimes where target and source servers dominate are not contiguous. Additionally, the phase transitions and the order in which they appear depend delicately on how distributed the data is and the quality of the source data. For a detailed account of all the phase transitions and their intricate dependence on the problem parameters, see the discussion following Corollary \ref{cor:pure-homogeneous-minimax}.
	
	In practical applications, we strive to borrow information from multiple sources which are very heterogeneous with respect to both the classification transfer quality, in terms of $\gamma$, and the privacy constraints. Take for example, the heart disease dataset of \cite{detrano1989international}, we want to predict the propensity of heart disease in a patient, based on their demographic information and clinical measurements. The data is collected from four different hospitals in Cleveland, Hungary, Long Beach, and Switzerland. The number of patients and disease prevalence are widely varied at each of these four hospitals. Moreover, a model fitted solely on one server shows drastically different performance on other servers, thus pointing to significant heterogeneity in the data. Details can be found in Section~\ref{sec:real-data}. At the same time, it is reasonable to assume each hospital having its own privacy tolerance level, determined by local guidelines. To tackle these challenges, in Section~\ref{sec:nonadaptive}, we develop minimax optimal methods in the heterogeneous setting, where the trade-offs between transfer learning and privacy are more nuanced. 
	
	In particular, we show in Section~\ref{sec:nonadaptive} that even when the source distributions are heterogeneous, with their own transfer parameters $\gamma_j$ and privacy constraints $\eps_j,\delta_j$ our weighted kernel based classifier continues to be minimax optimal. However the best choices of weights $w$ and bandwidth $h$ depends heavily on the knowledge of $\gamma_j$, which is typically unknown. The additional noise due to privacy compounds this issue even further, since choosing a worse bandwidth potentially implies adding a higher amount of noise than is required, thus unnecessarily worsening the performance. To alleviate this issue, we take a data-adaptive approach to choose the best bandwidth from a grid of possible choices. This is based on the popular Lepski method fine-tuned to our setting to accommodate the additional noise for privacy. Given the privacy parameters $(\eps,\delta)$, our adaptation algorithm outputs a classifier that attains the minimax rate, modulo a $O({\rm polylog}(n_*))$ factor, where $n_*$ is the sample size of the data pooled across all servers. Details can be found in Section~\ref{sec:adaptation}.
	
	The rest of the paper is organized as follows. In Section~\ref{sec:problem-form}, we provide the background and formulate our problem in detail. Section~\ref{sec:main-result}  presents the minimax rate of excess risk for our problem across various specific cases, as well as the most general case. Section~\ref{sec:nonadaptive} introduces our kernel based classifier, derives its excess risk bounds, and states the minimax lower bound. Section~\ref{sec:adaptation} describes the data-driven adaptive procedure for bandwidth and weight selection. We evaluate our proposed method and compare it to existing work via several numerical experiments on simulated and real datasets in Section~\ref{sec:numerical-studies}. The paper concludes with a discussion on possible future work in Section~\ref{sec:discussion}. All proofs can be found in the supplementary material~\cite{auddy2024supp}.
	
	%%%%%%%%%%%%%%%%%%%%%%%%%%%%%%%%%%%%%%%%%%
	\section{Problem Formulation}\label{sec:problem-form}
	%%%%%%%%%%%%%%%%%%%%%%%%%%%%%%%%%%%%%%%%%%
	
	In this section, we outline the general framework for transfer learning under distributed privacy constraints. Our dataset is distributed across \(m + 1\) servers, indexed by the set \(\{0,1,\ldots,m\}\). The dataset is categorized as target and source. On server \(0\) (also called the target server), we have \(n_0\) i.i.d. samples from the distribution \(P_0\), while on server \(j\) (the source servers) for \(j \in 1,\ldots,m\), we have \(n_{j}\) i.i.d. samples from distribution \(P_j\). All of the probability measures $\{P_j\}_{j=0}^m$ are defined on the measurable space $(\cZ,\scrZ)$. Let \(Z^{(0)} = \{Z_i^{(0)}\}_{i=1}^{n_{0}}\) denote the \(n_0\) realizations from \(P_0\) on the target server. Let us denote by \(Z^{(j)} = \{Z_i^{(j)}\}_{i=1}^{n_{j}}\) the \(n_{j}\) realizations from \(P_j\) on the \(j\)th source server for $j=1,\ldots,m$. These servers serve as the source data, and our goal is to learn the model for our target distribution \(P_0\).

	For each source server i.e $j=1,\ldots,m$, we send a (randomized) transcript \(\tilde{T}^{(j)}\) based on \(Z^{(j)}\) to the target server $0$, where the law of the transcript is given by a distribution conditional on \(Z^{(j)}\), \(\mathbb{P}(\cdot|Z^{(j)})\), on a measurable space \((\cT,\scrT)\). For $j = 1,\ldots,m$ the transcript \(\tilde{T}^{(j)}\) has to satisfy a \((\varepsilon_{j} , \delta_{j} )\)-differential privacy constraint.
	
	\begin{defn}\label{def:diff-privacy}
		The transcript \( \tilde{T}^{(j)} \) is \((\varepsilon_j, \delta_j)\)-differentially private if for all \( A \in \cA \) and \( z,z' \) differing in one individual datum, it holds that
		\[
		\mathbb{P}\left( \tilde{T}^{(j)} \in A \,|\, Z^{(j)} = z \right) \leq e^{\varepsilon_j} \mathbb{P}\left( \tilde{T}^{(j)} \in A \,|\, Z^{(j)} = z' \right) + \delta_j.
		\]
	\end{defn} 
	The target server can look at the private transcripts $\{\tilde{T}^{(j)}\}_{j=1}^m$ and the target data $Z^{(0)}$ while constructing the final private transcript $\tilde{T}$. Hence $\tilde{T}$ satisfies \((\varepsilon_{0} , \delta_{0} )\)-interactive differential privacy constraint, which is defined as follows:
	\begin{defn}\label{def:target-private}
		The transcript \( \tilde{T} \) is \((\varepsilon_0, \delta_0)\)-differentially private if for all \( A \in \cA \) and \( z,z' \) differing in one individual datum and for all $t_j \in \cT$ for $j=1,\ldots,m$, it holds that
		\begin{align*}
		&~\mathbb{P}\left( \tilde{T} \in A \,|\, Z^{(0 )} = z, \tilde{T}^{(j)} =t_j \text{ for }1\le j\le m\right)\\
		\leq&~ e^{\varepsilon_0} \mathbb{P}\left( \tilde{T} \in A \,|\, Z^{(0)} = z',  \tilde{T}^{(j)} =t_j \text{ for }1\le j\le m \right) + \delta_0.
		\end{align*}
	\end{defn}
	This privacy constraint can be understood as follows: if we condition on the outcome of all other servers then the distribution of the final private transcript $\tilde{T}$  does not change much if one of the datum on the target server is changed.
	
	In transfer learning, the focus is on scenarios where multiple parties, such as hospitals, possess heterogeneous data with differing underlying distributions. Employing distributed protocols in such contexts ensures differential privacy while yielding outputs from each participating party. Within this framework, transcripts generated by each source server rely solely on its local data, with no exchange of information occurring between source servers. Communication is solely between the source and target servers. Each of the source server transmits its transcripts to the target server. The target server utilizing all the transcripts $( \tilde{T}^{(1)},\ldots,\tilde{T}^{(m)})$ from the other servers and target data $Z^{(0)}$, computes the final private transcript $\tilde{T}$.  This scenario often arises when multiple trials involving a population similar to that of the target server are conducted, yet individual locations, such as hospitals, opt against consolidating their original data due to privacy apprehensions.
	
	In the context of transfer learning for nonparametric classification our data looks like a couple $Z_i^{(j)}:=(X_i^{(j)},Y_i^{(j)})$, for $i=1,\dots,n_{j}$; $j=1,\dots,m$ for the source servers, and $Z_i^{(0)}:=(X_i^{(0)},Y_i^{(0)})$, for $i=1,\dots,n_{0}$ for the target server. We assume that $Z_i^{(j)}$ takes values in $\calZ:=[0,1]^d\times\{0,1\}$. We regard $X\in [0,1]^d$ as a vector of features corresponding to an object and $Y\in\{0,1\}$ as a label indicating that the object belongs to one of two classes. Our goal is to propose distributed DP protocols $\tilde{T}^{(j)}$ for each server and construct classifier $\widehat f:[0,1]^d \to \{0,1\}$ based on the final private transcript $\{\tilde{T}\}$. Unlike the traditional federated learning framework, there's no central server; alternatively, we can consider the target server as acting in a central capacity. 
	We denote the vector of privacy budgets as $(\bm\varepsilon,\bm\delta) =\{(\varepsilon_j,\delta_j)\}_{j=0}^m$ and the class of distributed DP classifiers $\widehat f$ by $\cM_{\bm \varepsilon,\bm\delta}$. Next we denote
	\begin{align*}
		\eta_{j}(X^{(j)}):=&~ \PP(Y^{(j)} = 1|X^{(j)})\text{ for the source servers }j=1,\dots,m;\text{ and } \\
		\eta_{0}(X^{(0)}):=&~ \PP(Y^{(0)} = 1|X^{(0)})\text{ for the target server},
	\end{align*}
	as the (source and target) regression functions of $Y$ on $X$. We denote the marginal distribution of $X$ for the $j$th server, $j=0,\ldots,m$ as $P_{j}^X$. Define the classification error of a classifier $f$ under the target distribution $P_0$ as 
	\[
	R_0( f) := P_0(Y \neq f(X))
	\]
	The Bayes decision rule is a minimizer of the of the risk $R_0(f)$ which has the form $f^*_0(X) = \bbm{1}\{\eta_0(X) \geq 1/2\}$.
	The goal of transfer learning is to transfer the knowledge gained from the source data together with the information in the target data to construct a classifier which minimizes the excess risk on the target data
	\[
	\calE_0(\hat{f})
	=\EE [R_0(\hat{f})]-R_0(f^*_0)
	\]
	%=\EE(|2\eta_Q(X)-1|\bbm{1}_{\hat{f}(X)\neq f^*(X)}),
	
	Under the posterior drift model we quantify the similarity between the regression functions $\{\eta_j\}_{j=1}^m$ and $\eta_0$ as follows:

	\begin{defn}[\textbf{Relative Signal Exponent (RSE)}] The class $\Gamma(\bm\gamma,C_{\bm \gamma})$ with relative signal exponent $\bm \gamma = (\gamma_1,\ldots,\gamma_m) \in \RR^m_+$ and constants $C_{\bm \gamma} =  (C_1,\ldots,C_m) \in \RR^m_+$, is the set of distribution tuples $(P_0,P_1,\ldots,P_m)$ that satisfy for $1\le j\leq m$  
		\begin{enumerate}
			\item ${\rm sign}\left(\eta_{j}(x)-\frac{1}{2}\right)={\rm sign}\left(\eta_0(x)-\frac{1}{2}\right)$ for all $1\le j\le m$ and all $x\in[0,1]^d$.
			\item $\left\vert\eta_{j}(x)-\frac{1}{2}\right\vert\ge 
			C_j
			\left\vert
			\eta_{0}(x)-\frac{1}{2}
			\right\vert^{\gamma_j}$ for some $\gamma_j>0$, for all $1\le j\le m$ and all $x\in[0,1]^d$.
		\end{enumerate}
	\end{defn}
	\begin{remark}{\rm
			The first part follows from the assumption that the Bayes classifier $f^*$ is the same for both source and target populations. The second part introduces a parameter $\gamma$ which controls the signal strength of the source data from $P$, with respect to the target data from $Q$. See Definition 1 and Remark 1 of \cite{cai2021transfer}.
	}\end{remark}
	In addition to the RSE assumption we also need to assume smoothness of $\eta_0$ and characterize it behavior  near $1/2$.
	These assumptions are standard in nonparametric classification and was first introduced in \cite{audibert2007plugin}. 
	\begin{defn}[\textbf{H{\"o}lder Smoothness}] The regression function $\eta_0$ belongs to the H{\"o}lder class of functions denoted by $\Sigma(\beta,L)$ $(0<\beta \le 1)$ which is defined as the set of functions satisfying:
		$$
		|\eta_0(x)-\eta_0(x')|\le L\|x-x'\|^{\beta}
		\quad
		\text{for } x,x'\in [0,1]^d.
		$$
	\end{defn}
	\begin{defn}[\textbf{Margin Assumption (MA)}] \label{def:margin-assn} The margin class $\cM(\alpha,C_\alpha)$ with $\alpha \ge 0$ and $C_\alpha >0$ is defined as the set of distributions $P_0$ such that 
		$$
		P_0^X(0\le |\eta_0(X)-1/2|\le t)
		\le C_\alpha t^{\alpha}
		\,\,
		\text{for all }
		t>0.
		$$
	\end{defn}
	
%	\begin{remark}\label{rem:beta<1}{\rm
%			The smoothness and the margin assumptions are standard in nonparametric classification and was first introduced in \cite{audibert2007plugin}. 
%			%We will focus our attention in this paper to the case $0<\beta\leq 1$. The discussion of the case where the smoothness parameter is greater than $1$ is deferred to the discussion section. 
%	}\end{remark}
	
	Another definition is about marginal density of $X$, $P_j^X$ for $j=1,\ldots,m$.
	\begin{defn}[\textbf{Common Support and Strong Density Assumption (SD)}] We assume that $P_j^X$ for $j=0,\ldots,m$  have the identical support on a compact $(c_\mu, r_\mu)$ regular set $A\subset [0,1]^d$ and has a density $g_j$ w.r.t. the Lebesgue measure bounded away from zero and infinity on $A$:
		$$
		g_{\min} \le g_j(x) \le g_{\max}\text{ for }x\in A \text{ and }g_j(x) =0\text{ otherwise},
		$$
		where $c_0,r_0>0$ and $0< g_{\min} < g_{\max} <\infty$ are fixed constants.
		We denote the set of marginal distributions $(P_0^X,\ldots,P_m^X)$ which satisfy the above constraints as $\cS(\mu,c_\mu,r_\mu)$ where $\mu = (g_{\min},g_{\max})$.
	\end{defn}	
	
	\begin{remark}{\rm
			In this paper we focus our attention to the case when the marginal densities have regular support and are bounded from below and above on their support. Moreover we assume that $\alpha \beta \leq d$ throughout the  paper. This is because in the other regime $(\alpha\beta > d)$, there is no distribution $P_0^X$ such that the regression function $\eta_0$ crosses $1/2$ in the interior of the support of $P_0^X$ (\cite{audibert2007plugin}) and hence  this case only contains the trivial cases for classification.  
	}\end{remark}
	
	We put all the definitions together to define the class of distributions we consider in the posterior drift model as
	\begin{align*}
		&~\Pi(\bm \gamma,C_{\bm \gamma}, \beta, L, \alpha, C_\alpha, \mu, c_\mu,r_\mu) \\
		:=&~ \{ (P_0, P_1, \ldots,P_m) : (P_0, P_1, \ldots,P_m) \in \Gamma(\bm \gamma,C_{\bm \gamma}), 
		\eta_0 \in \Sigma(\beta,L),\\
		&\hspace{120pt} \P_0^X \in \cM(\alpha,C_\alpha), (P_0^X, P_1^X, \ldots,P_m^X) \in  \cS(\mu,c_\mu,r_\mu)\}
	\end{align*}
	For the rest of the paper we will use the shorthand $\Pi(\alpha,\beta,\bm \gamma,\mu)$ or $\Pi$ if there is no confusion.
	
	%%%%%%%%%%%%%%%%%%%%%%%%%%%%%%%%%%%%%%%%%%
	\section{Main Results}\label{sec:main-result}
	%%%%%%%%%%%%%%%%%%%%%%%%%%%%%%%%%%%%%%%%%%
	
	In this section, we present the key findings of our paper, where we establish the minimax rate of convergence for  transfer learning under differential privacy constraints, specifically addressing the nonparametric classification problem.We divide our results into two subsections: Section \ref{subsec:homogeneity} covers the homogeneous case, while Section  \ref{subsec:general} addresses the general heterogeneous case.
	
	%%%%%%%%%%%%%%%%%%%%%%%%%%%%%%%%%%%%%%%%%%
	\subsection{Minimax Rates under Source Homogeneity}\label{subsec:homogeneity}
	%%%%%%%%%%%%%%%%%%%%%%%%%%%%%%%%%%%%%%%%%%
	
	To derive meaningful and interpretable insights from our minimax rate, we first examine the scenario where the source servers are exchangeable in terms of the distributed classification problem under transfer learning and privacy constraints. This homogeneous scenario is characterized by equal sample sizes ($n_j = n$), privacy parameters ($\varepsilon_j = \varepsilon$, $\delta_j = \delta$) and transfer exponents ($\gamma_j = \gamma$) for all $j=1,\ldots,m$.
	
	\begin{theorem}\label{th:homogeneous-minimax-rate}
		Suppose $n_j =n,\varepsilon_j =\varepsilon,\delta_j = \delta$ and $\gamma_j =\gamma$ for all $j=1,\ldots,m$ and assume that $\delta =o((nm)^{-1})$. Then the minimax rate for the excess risk satisfies
		\begin{align*}
			\inf_{\widehat{f}\in \cM(\bm\varepsilon,\bm\delta)}\sup_{(P_0,\ldots,P_m) \in \Pi}  \cE_0(\widehat f) 
			\asymp\bigg[ L_N\bigg\{
			&~
			\left( n_0^{\frac{1}{2\beta + d}} \wedge (n_0^2\varepsilon^2_0)^{\frac{1}{2\beta + 2d}}\right) \\
			&~+ \left( (mn)^{\frac{1}{2\beta \gamma+ d}} 
			\wedge (mn^2\varepsilon^2)^{\frac{1}{2\beta \gamma+ 2d}}\right) \bigg\}^{-\beta(1+\alpha)}
			\wedge 1\bigg]
		\end{align*}
		for a sequence $L_N$ of order at most $(\log((\delta\wedge \delta_0)^{-1}))^{\frac{\beta(1+\alpha)}{2\beta(\gamma\wedge 1)+d}}$.
	\end{theorem}
	
	Note that the minimax rate depends on the sum of quantities: the first of which determines the minimax rate for the problem using only the target data, while the second corresponds to the minimax rate for the problem using solely the source data from the source servers.  Some remarks are in order regarding the minimax rate obtained of Theorem~\ref{th:homogeneous-minimax-rate} in some special settings:
	\begin{remark}[Private classification with no transfer]{\rm
		When the number of servers $m=0$ or there is no source data, i.e., $n=0$ for all source servers we obtain the single server minimax classification rate under the $(\eps, \delta)$-DP constraint, given by 
		\[
		\inf_{\widehat{f}\in \cM(\bm\varepsilon,\bm\delta)}\sup_{(P_0,\ldots,P_m) \in \Pi}  \cE_0(\widehat f) 
		\asymp
		\bigg[ L_N
		\left( n_0^{\frac{1}{2\beta + d}} \wedge (n_0^2\varepsilon^2_0)^{\frac{1}{2\beta + 2d}}\right)^{-\beta(1+\alpha)}
		\wedge 1
		\bigg].
		\]
	}\end{remark}
	
	\begin{remark}[Non-private transfer learning]{\rm
		When the privacy requirements are not stringent, the tradeoff is completely characterized by a comparison between the non-private rates of the target and the source. In particular, if $\eps>\left(m^{d/2}n^{-\beta\gamma}\right)^{\frac{1}{2\beta\gamma+d}}$, we find the non-private transfer learning rates
		\[
		\inf_{\widehat{f}\in \cM(\bm\varepsilon,\bm\delta)}\sup_{(P_0,\ldots,P_m) \in \Pi}  \cE_0(\widehat f) 
		\asymp
		\left(n_0+(mn)^{\frac{2\beta+d}{2\beta\gamma+d}}\right)^{-\frac{\beta(1+\alpha)}{2\beta+d}}
		\]
		which coincides with the results of %Theorems 1 and 2 in
		\cite{cai2021transfer}, in the transfer homogeneous regime with equal source sample sizes. 
	}\end{remark}

	In the current setting, we further emphasize how the requirement of privacy leads to a worsening of the rate if the same amount of data is distributed across a larger number of servers. Since the sources are all equivalent in terms of data quality (as quantified by $\gamma$), traditional knowledge suggests a rate depending on the pooled source sample size $mn$. In the $m$-server non-private transfer learning setup of \cite{cai2021transfer}, that is indeed the case, as also shown by the non-private rate of $(mn)^{-\frac{\beta(1+\alpha)}{2\beta\gamma+2d}}$ appearing in Theorem~\ref{th:homogeneous-minimax-rate}. Note however that the private source rate is given by $m^{\frac{\beta(1+\alpha)}{2\beta\gamma+2d}}(mn\eps)^{-\frac{\beta(1+\alpha)}{\beta\gamma+d}}$, so that if the pooled sample size $mn$ remains fixed, the rate worsens as $m$, the number of servers, increases. This phenomenon is reminiscent of the minimax rate behavior observed in \cite{Cai2023FL-NP-Regression} in the non-private setting. 
	
	In order to further our understanding about interplay between the transfer exponent $\gamma$ and privacy parameters we restrict our attention to the case where the target server has same privacy budget, i.e., $\varepsilon_{0} =\varepsilon,\,\delta_0 = \delta$ and the number of target samples is between $n$ and $mn$. Other sample size regimes can be described similarly. As is clear from Theorem~\ref{th:homogeneous-minimax-rate}, the minimax rate of decay for the excess risk is given in this case by:
	\begin{equation}\label{eq:minimax-rate-homogenous}
		\cE_0(\widehat f) \asymp\left[ L_N\left\{\left( n_0^{-\frac{\beta(1+\alpha)}{2\beta + d}} \vee (n_0^2\varepsilon^2)^{-\frac{\beta(1+\alpha)}{2\beta + 2d}}\right) \bigwedge \left( (mn)^{-\frac{\beta(1+\alpha)}{2\beta \gamma+ d}} \vee (mn^2\varepsilon^2)^{-\frac{\beta(1+\alpha)}{2\beta \gamma+ 2d}}\right) \right\}\wedge 1\right].
	\end{equation}
	Figure~\ref{fig:theo-rates} illustrates the different regimes showing which of the four terms on the right hand side of \eqref{eq:minimax-rate-homogenous} ends up determining the overall rate of excess risk.

	We will refer to the four terms on the right hand side in \eqref{eq:minimax-rate-homogenous} as the non-private target rate (${\rm NP}_t$), the private target rate (${\rm P}_t$), the non-private source rate (${\rm NP}_s$), and the private source rate ($P_s$) respectively. Depending on the value of the common privacy parameter $\eps$ and the transfer exponent $\gamma$, the overall rate will be determined by one of these four rates, as demonstrated by the following table. Corollary~\ref{cor:pure-homogeneous-minimax} formally states the results of Table~\ref{tab:rates} along with the endpoints of the various transfer and privacy regimes.

	\begin{table}
		\caption{\label{tab:rates} Minimax rate of excess risk at different transfer and privacy regimes. See Corollary~\ref{cor:pure-homogeneous-minimax}.}
		%\centering
		%\caption{\label{tab:rates} Table1}
		%\caption{\label{tab:rates} Minimax rate of excess risk at different transfer and privacy regimes. Here we write ${\rm P}_s=L_N(mn^2\eps^2)^{-\frac{\beta(1+\alpha)}{2\beta\gamma+2d}}$, ${\rm P}_t=L_N(n_0^2\eps^2)^{-\frac{\beta(1+\alpha)}{2\beta+2d}}$, ${\rm NP}_s=L_N(mn)^{-\frac{\beta(1+\alpha)}{2\beta\gamma+d}}$, ${\rm NP}_t=L_Nn_0^{-\frac{\beta(1+\alpha)}{2\beta+d}}$ to denote the four kinds of possible rates. The endpoints of the different regimes can be found in Corollary~\ref{cor:pure-homogeneous-minimax}.}
		%\centering
		%\fbox{
			%\hspace{-7cm}
			\begin{tabular}{|c|c|cl|clc|}
			\hline
			\multirow{2}{*}{\diagbox{\small Transfer}{\small Privacy}} & \multirow{2}{*}{$\eps\in (0,\eps^{(1)}]$}  & \multicolumn{2}{c|}{\multirow{2}{*}{$\eps\in (\eps^{(1)},\eps^{(2)}]$}}   & \multicolumn{2}{c|}{\multirow{2}{*}{$\eps\in (\eps^{(2)},\eps^{(3)}]$}}    & \multirow{2}{*}{$\eps\in (\eps^{(3)},1]$}    \\
			&                    & \multicolumn{2}{c|}{}                    & \multicolumn{2}{c|}{}                     &                      \\ \hline
			\multirow{3}{*}{$\gamma\in (0,1]$} &   \multirow{8}{*}{1}                 & \multicolumn{2}{c|}{\multirow{3}{*}{
					$\begin{cases}
						{\rm P}_t \,\,&\text{ if }\eps\le \eps^{(11)}\\
						{\rm P}_s \,\,&\text{ if }\eps>\eps^{(11)}
					\end{cases}
					$
			}} & \multicolumn{3}{c|}{\multirow{3}{*}{${\rm NP}_s$}}                        \\ 
			&                    & \multicolumn{2}{c|}{}                    & \multicolumn{3}{c|}{}                                            \\
			&                    & \multicolumn{2}{c|}{}                    & \multicolumn{3}{c|}{}                                            \\ \cline{1-1} \cline{3-7}  
			\multirow{3}{*}{$\gamma\in (1,\gamma^{(*)}]$} &                    & \multicolumn{2}{c|}{\multirow{5}{*}{
					$\begin{cases}
						{\rm P}_s \,\,&\text{ if }\eps\le \eps^{(11)}\\
						{\rm P}_t \,\,&\text{ if }\eps>\eps^{(11)}
					\end{cases}
					$
			}}                    & \multicolumn{2}{c|}{\multirow{3}{*}{$
					\begin{cases}
						{\rm NP}_t \,\,&\text{ if }\eps\le \eps^{(21)}\\
						{\rm P}_s \,\,&\text{ if }\eps>\eps^{(21)}
					\end{cases}
					$}} & \multirow{3}{*}{${\rm NP}_s$} \\
			&                    & \multicolumn{2}{c|}{}                    & \multicolumn{2}{c|}{}                     &                      \\ 
			&                    & \multicolumn{2}{c|}{}                    & \multicolumn{2}{c|}{}                     &                      \\  \cline{1-1} \cline{5-7}  
			\multirow{2}{*}{$\gamma\in (\gamma^{(*)},\infty)$} &  & \multicolumn{2}{c|}{} & \multicolumn{3}{c|}{\multirow{2}{*}{${\rm NP}_t$}}                        \\    
			&                    & \multicolumn{2}{c|}{}                    & \multicolumn{3}{c|}{}\\              	
			\hline
		\end{tabular}
	\end{table}
	
	\begin{corollary}\label{cor:pure-homogeneous-minimax}
		Suppose $n_j = n, \gamma_j =\gamma \,\forall \, 1\leq j\leq m$,  $n\le n_0\le mn$, and equal privacy budget $\varepsilon_j=\varepsilon,\,\delta_j = \delta  \,\forall \, 0\leq j\leq m$ . Further assume that $\delta =o((mn)^{-1})$. Then the minimax rate for the excess risk are as given in Table~\ref{tab:rates} with the various regimes characterized by the following endpoints:
		\begin{enumerate}
			\item $\gamma^{(*)}=\frac{1}{2\beta}\left[\frac{(2\beta+d)\log mn}{\log n_0}-d\right]$.
			\item $\eps^{(1)}=(\sqrt{m}n)^{-1}\wedge n_0^{-1}$; $\eps^{(2)}=n_0^{-\frac{\beta}{2\beta+d}}$; $\eps^{(3)}=\left(m^{d/2}n^{-\beta\gamma}\right)^{\frac{1}{2\beta\gamma+d}}$.
			\item $\eps^{(11)}=
			\begin{cases}
				\left[(\sqrt{m}n)^{\beta+d}n_0^{-(\beta\gamma+d)}\right]^{\frac{1}{\beta(\gamma-1)}}%\vee \eps^{(1)}
				\quad&\text{if }\gamma\neq 1,\\
				\eps^{(1)}\quad&\text{if }\gamma= 1,\,n_0\le\sqrt{mn^2},\\
				\eps^{(2)}\quad&\text{if }\gamma= 1,\,n_0>\sqrt{mn^2}.
			\end{cases}
			$
			\item $\eps^{(21)}=(\sqrt{m}n)^{-1}n_0^{\frac{\beta\gamma+d}{\beta+d}}$.
		\end{enumerate}
	\end{corollary}

	As we increase the value of $\eps$, we observe interesting phenomena characterized by distinct phase transitions. Not surprisingly, the rates behave differently based on whether $\gamma$ is small --- where the quality of the source data is relatively better than the target data --- versus when $\gamma$ is larger. The two transition points in $\gamma$ are at $\gamma=1$ and at  $\gamma^{(*)}=\frac{1}{2\beta}\left[\frac{(2\beta+d)\log mn}{\log n_0}-d\right]$. The different rates can be described based on ultra-high, high, moderate, and low privacy regimes.
	
	First, in the ultra-high privacy regime where $0<\eps\le\eps^{(1)}$, the privacy requirements are so severe that no classifier has disappearing excess risk in this regime, and a random guess is the best one can do. Next, we move to the high privacy regime of $\eps^{(1)}<\eps\le \eps^{(2)}$, where the private rates of both the target and the source dominate over their non-private counterparts. Another interesting phenomenon emerges based on whether $\gamma$ is smaller than or greater than one. If $\gamma<1$ and the privacy requirement is high, the target private rates appear for very small $\eps$, followed by the source private rates. This is because when $\gamma<1$, the source data are of better quality and hence for sufficiently small $\eps$, the noise added to the target dominates over the noise for the sources. The pattern reverses when $\gamma>1$ where the source private rates appear first. 
	
	Further increasing $\eps$, we arrive at the moderate privacy regime where $\eps^{(2)}<\eps\le\eps^{(3)}$. For small $\gamma$, characterized by $\gamma\le 1$, the high quality source data points prove to be particularly beneficial. This results in the non-private source rate dominating over the non-private target rate, as well as the noises added for privatizing the source and the target. In contrast, when $\gamma$ is very large, in particular $\gamma>\gamma^{(*)}$, the source data are particularly poor and not useful, so that the target non-private rate therefore dominates over all the other contenders. A more interesting picture emerges for moderate $\gamma$ given by $1<\gamma\le \gamma^{(*)}$. In this regime of intermediate transfer and medial privacy, we find the non-private target rate first, followed by private source rates. 
	
	Finally for $\eps>\eps^{(3)}$, the low privacy scenario emerges. Here the effect of extraneous noise for privacy is not at all significant. We therefore only find the non-private rates based on the transfer efficacy. For $\gamma\le \gamma^{(*)}$, the source data are more useful, resulting in ${\rm NP}_s$ governing the overall rate. As expected, for $\gamma>\gamma^{(*)}$ the relatively poorer source data do not contribute to the transfer. Consequently the non-private rates from the target becomes dominant, reflecting the diminished relevance of the source data in such scenarios.

	In the rest of this subsection, we describe another specialized setting where we allow one of the source servers to be public. To demonstrate the effect of publicly available data, we take $m=2$ sources, with one private and one public source server. The minimax rate for excess risk is then given by the following corollary:
	
	\begin{corollary}\label{cor:public-sources}
		Suppose that $\gamma_1=\gamma_2=\gamma$, $\eps_1=\infty$, $\eps_0=\eps_2=\eps$ and $n_2>n_0^{\frac{2\beta\gamma+d}{2\beta+d}}$, $\delta_0\vee \delta_2 =o(n_2^{-1})$. Then the minimax rate for the excess risk satisfies the following.
		\begin{enumerate}
			\item If $n_1>n_2$, 
			$\inf_{\widehat{f}\in \cM(\bm\varepsilon,\bm\delta)}\sup_{(P_0,\ldots,P_m) \in \Pi}  \cE_0(\widehat f) 
			\asymp
			n_1^{-\frac{\beta(1+\alpha)}{2\beta\gamma+d}}
			$.
			\item If $n_0^{\frac{2\beta\gamma+d}{2\beta+d}}\le n_1\le n_2$, then
			\begin{align*}
				\inf_{\widehat{f}\in \cM(\bm\varepsilon,\bm\delta)}\sup_{(P_0,\ldots,P_m) \in \Pi}  \cE_0(\widehat f) 
				\asymp
				\begin{cases}
					L_Nn_1^{-\frac{\beta(1+\alpha)}{2\beta\gamma+d}} \qquad
					&\text{if }\eps\le n_2^{-1}n_1^{\frac{\beta\gamma+d}{2\beta\gamma+d}}\\
					L_N(n_2^2\eps^2)^{-\frac{\beta(1+\alpha)}{2\beta\gamma+2d}} \qquad
					&\text{if } n_2^{-1}n_1^{\frac{\beta\gamma+d}{2\beta\gamma+d}} < \eps \le 
					L_Nn_2^{-\frac{\beta\gamma}{2\beta\gamma+d}}
					\\
					L_Nn_2^{-\frac{\beta(1+\alpha)}{2\beta\gamma+d}} \qquad
					&\text{if } n_2^{-\frac{\beta\gamma+d}{2\beta\gamma+d}} < \eps \le 
					1.
				\end{cases}
			\end{align*}
			\item If $n_1 \le n_0^{\frac{2\beta\gamma+d}{2\beta+d}}\le n_2$, then
		\end{enumerate}
		\begin{align*}
			&~\inf_{\widehat{f}\in \cM(\bm\varepsilon,\bm\delta)}\sup_{(P_0,\ldots,P_m) \in \Pi}  \cE_0(\widehat f)\\ 
			\asymp &~
			\begin{cases}
				L_Nn_1^{-\frac{\beta(1+\alpha)}{2\beta\gamma+d}} 
				\qquad
				&					\text{if }\eps\le \tilde{n}^{-\beta}\\
				\bigg[L_N\bigg\{			
				\left( n_0^{\frac{1}{2\beta + d}} \wedge (n_0^2\varepsilon^2)^{\frac{1}{2\beta + 2d}}\right)
				+ \left( n_2^{\frac{1}{2\beta \gamma+ d}} 
				\wedge (n_2^2\varepsilon^2)^{\frac{1}{2\beta \gamma+ 2d}}\right) \bigg\}^{-\beta(1+\alpha)}
				\wedge 1\bigg]
				&\text{otherwise,}
			\end{cases}
		\end{align*}
		where $\tilde{n}=n_0^{\frac{1}{2\beta+d}}\wedge 
		n_2^{\frac{\gamma}{2\beta\gamma+d}}$.	Here $L_N$ is a sequence of order at most $(\log((\delta_0\wedge \delta_2)^{-1}))^{\frac{\beta(1+\alpha)}{2\beta(\gamma\wedge 1)+d}}$.
	\end{corollary}
	
	The above corollary demonstrates three different rates based on $n_1$, the size of the public source server, with respect to the private source and the target. When $n_1$ is larger than $n_2$, the public source has a significantly large number of samples, enough for the non-private rates from this server to dominate over the rest. Next, when $n_1$ is relatively large with respect to the target sample size, but still smaller than the private source server, we find three different rates. In the high privacy regime of very small $\eps$, the availability of public data ensures that we do not have to pay an additional price of private rates by adding noise. Next, the private rates of the source appear: the reason being that the private source has higher number of samples, and in this intermediate privacy regime, even after adding additional noise, the private source turns out to be more useful than the public server. For even lower privacy, the extraneous noise level diminishes further and the non-private rate from source server 2 appears. The final case is when the public data size $n_1$ is relatively smaller than the sample size in the two privacy constrained servers. In the regime of very high privacy, the public data is still effective, and enables one to avoid the private rates. For moderate and low privacy, the tradeoff is characterized by the relative signal strengths and privacy requirements of the private source and target servers. The exact rates would be similar to Table~\ref{tab:rates} in this scenario. 
	
	%%%%%%%%%%%%%%%%%%%%%%%%%%%%%%%%%%%%%%%%%%
	\subsection{Minimax Rates in General Setting}\label{subsec:general}
	%%%%%%%%%%%%%%%%%%%%%%%%%%%%%%%%%%%%%%%%%%
	
	We now turn our attention to the general case where the sample sizes $n_j$, transfer exponents $\gamma_j$, privacy parameters $(\eps_j,\delta_j)$ are all allowed to vary for $0\le j\le m$. Our main result, captured in Theorem \ref{thm:master-minimax-rate}, quantifies the rate. The homogeneous case described earlier can be thought of as a special case of this vastly more general setting.
	
	\begin{theorem}\label{thm:master-minimax-rate}
		Let $r \in \RR_+$ be the solution to the following equation:
		
		\begin{equation}\label{eq:optimal-bandwidth-choice}
			(n_0\wedge n_0^2\eps_{0}^2r^d)r^{2\beta+d}
			+
			\displaystyle\sum_{j=1}^m
			(n_{j}\wedge n_{j}^2\eps_{j}^2r^d)
			r^{2\beta\gamma_j+d}
			=~1
		\end{equation}
		The minimax rate for excess risk is given by 
		\begin{equation}
			\inf_{\widehat{f}\in \cM(\bm\varepsilon,\bm\delta)}\sup_{(P_0,\ldots,P_m) \in \Pi(\alpha,\beta,\bm \gamma,\mu)}  \cE_0(\widehat f)\asymp\left( L_N r^{\beta(1+\alpha)} \wedge 1\right).
		\end{equation}
		whenever $\sum_j n_j\delta_j \to 0$, for a sequence $L_N$ of order at most $\left(
		-
		\log\left(
		\delta_{\min}
		\right)
		\right)^{\frac{\beta(1+\alpha)}{2\beta\gamma_{\min}+d}}$.
	\end{theorem}
	It is important to note that ~\eqref{eq:optimal-bandwidth-choice} always yields a positive solution, since the left-hand side of \eqref{eq:optimal-bandwidth-choice} is a strictly increasing continuous function of $r$, ranging from 0 to $\infty$ as $r$ varies from $0$ to $\infty$. Thus the function of $r$ defined on the left side of \eqref{eq:optimal-bandwidth-choice} must take the value $1$ somewhere in $\mathbb{R}_+$. We now provide a brief commentary on the derived result, starting with a comparison with non-private transfer learning.
	
	\begin{remark}[General non-private transfer learning]{\rm
		When the privacy budget is large, setting $\varepsilon_j =\infty$ for $j=1,\ldots,m$ in \eqref{eq:optimal-bandwidth-choice} we recover the non-private rate
		\[
		\inf_{\widehat{f}\in \cM(\bm\varepsilon,\bm\delta)}\sup_{(P_0,\ldots,P_m) \in \Pi(\alpha,\beta,\bm \gamma,\mu)}  \cE_0(\widehat f)\asymp
		\left(
		n_0+\sum_{j=1}^mn_j^{\frac{2\beta+d}{2\beta\gamma_j+d}}		
		\right)^{-\frac{\beta(1+\alpha)}{2\beta+d}}
		\]
		when $m$ is fixed. This coincides with Theorems 5 and 6 of \cite{cai2021transfer}. Note however that Theorem~\ref{thm:master-minimax-rate} allows the number of servers, $m$, to grow to $\infty$. 
	}\end{remark}

	In the most general case, the minimax optimal rate of convergence exhibited in the transfer learning problem under distributed privacy is determined by $r$, which has implicit dependencies on various factors, including the number of servers, privacy parameters, transfer exponents, and sample sizes. The value of $r$ determines the radius within which local methods can share information. In Section~\ref{sec:class-up-bd} we find that when using kernel estimators, the optimal bandwidth choice is innately connected to $r$, and in fact differs from the solution to \eqref{eq:optimal-bandwidth-choice} by at most a logarithmic factor. To further interpret the role of $r$, note that \eqref{eq:optimal-bandwidth-choice} quantifies the exact contribution of each server to the entire classification procedure. For each $j\in\{0,\dots,m\}$, the variance of the local estimator for the $j^{\text{th}}$ server is determined by the sum of two quantities: the inverse of the sample size, and the variance of the additional noise required for privacy. The term $n_j\wedge n_j^2\eps_j^2r^d$ then appears as a quantity proportional to the inverse of this variance and determines the \emph{precision} of the $j^{\text{th}}$ server. 
	
	%%%%%%%%%%%%%%%%%%%%%%%%%%%%%%%%%%%%%%%%%%
	\section{Minimax Optimal Classification Procedure}\label{sec:nonadaptive}
	%%%%%%%%%%%%%%%%%%%%%%%%%%%%%%%%%%%%%%%%%%
	
	In this section, we propose an optimal classifier and establish the minimax rate of convergence. In the first subsection, we develop a nonparametric classifier for the target population that appropriately utilizes information from the sources while satisfying privacy requirements for each server. In the second subsection, we prove a minimax lower bound, demonstrating that our classifier achieves the minimax rate of convergence in the distributed private transfer learning context.
	
	%%%%%%%%%%%%%%%%%%%%%%%%%%%%%%%%%%%%%%%%%%
	\subsection{Classifier}\label{sec:class-up-bd}
	%%%%%%%%%%%%%%%%%%%%%%%%%%%%%%%%%%%%%%%%%%
	
	We now describe a classifier for transfer learning with distributed privacy. Our method has three main steps. First, we use a kernel estimator to estimate $(\eta_j(x)-\frac{1}{2})g(x)$ for $j=0,1,\ldots,m$. Second, then use a convex combination of these estimators, where the weights are designed to borrow strength from the source servers under the transfer learning setup. The third step is to add a Gaussian noise to the weighted kernel estimator to satisfy privacy requirements. Our classifier is given by the sign of the noise perturbed weighted estimator.
	
	Consider a kernel $K(t)$ supported on $[-1,1]^d$ with the following properties:
	\begin{enumerate}
		\item $\int K(t)dt =1$.
		\item $K(\cdot)$ is $L_K$-Lipschitz.
		\item $\underset{t\in [-1,1]^d}{\max} K(t) \le c_K$.
		\item $\min\{K(t):\|t\|\le 1/2\} \ge b_K$.
		\item $K$ is positive definite.
	\end{enumerate}
	Let $X=x_0$ be the test point we wish to classify. Suppose first that $x_0\in[0,1]^d$. Then for the source samples we compute
	\begin{equation}\label{eq:def-Tnj-np}
		T^{(j)}_h(x_0)
		:= \dfrac{1}{n_{j}h^d}\sum_{i=1}^{n_{j}}
		\left(Y_i^{(j)}-\frac{1}{2}\right)
		K\left(
		\dfrac{X_i^{(j)}-x_0}{h}
		\right)
	\end{equation}
	for $0\le j\le m$. To satisfy privacy requirements we follow the framework of \cite{hall2013differential} and a Gaussian process to the above estimator. Note that $T^{(j)}_h$ belong to the RKHS $\calK$ given by linear combination $\{\sum_i\theta_iK((X_i-x)/h):\theta_i\in\RR\}$. For two functions $f=\sum_i\theta_iK((X_i-x)/h)$ and $g=\sum_i\tau_iK((X_i-x)/h)$, their inner product under $\calK$ is given by 
	$$
	\langle f,g\rangle_{\calK}=\sum_i\sum_j\theta_i\tau_jK\left(\dfrac{X_i-X_j}{h}\right).
	$$
	Let $T^{(j)\prime}_h$ be the versions of $T^{(j)}_h$ with $(X_1^{(j)},Y_1^{(j)})$ replaced by $(X_1^{(j)\prime},Y_1^{(j)\prime})$ for $j=0,1,\dots,m$. Then the RKHS norm of $T^{(j)}_h(\cdot)-T^{(j)\prime}_h(\cdot)$ can be bounded by 
	\begin{align*}\label{eq:Tj-rkhs-norm}
		~\Vert
		T^{(j)}_h-T^{(j)\prime}_h\Vert_{\calK}
		%	=&~\sqrt{
			%		\dfrac{1}{n_{P_j}^2h^{2d}}
			%		\left[
			%		\dfrac{1}{2}
			%		K(0)
			%		-
			%		2\left(Y_1^{(j)}-\frac{1}{2}\right)\left(Y_{1}^{(j)'}-\frac{1}{2}\right)
			%		K\left(\dfrac{X^{(j)}_{1}-X^{(j)'}_{1}}{h}\right)
			%		\right]}
		\le~ \dfrac{\sqrt{c_K}}{n_{P_j}h^d}
		%	\\
		\quad
		\text{and }
		\quad 
		~\Vert
		T^{(0)}_h-T^{(0)\prime}_h
		\Vert_{\calK}
		%	=&~\sqrt{
			%		\dfrac{1}{n_{Q}^2h^{2d}}
			%		\left[
			%		\dfrac{1}{2}
			%		K(0)
			%		-
			%		2\left(Y_{1}-\frac{1}{2}\right)\left(Y_{1}'-\frac{1}{2}\right)
			%		K\left(\dfrac{X_{1}-X_{1}'}{h}\right)
			%		\right]
			%	}
		\le~ \dfrac{\sqrt{c_K}}{n_{Q}h^d}
		.\numberthis
	\end{align*}
	Let us define $(m+1)$ independent mean zero Gaussian processes $\xi^{(j)}(\cdot)$ with covariance kernels 
	$$
	{\rm Cov}(\xi^{j}(s),\xi^{(j)}(t))=K\left(\frac{s-t}{h}\right)
	\text{ for }
	s,t\in[0,1].
	$$
	and
	\begin{equation}\label{eq:def-tild-xi}
		\tilde{\xi}_h{(\cdot)}
		=
		\dfrac{\sqrt{2c_K\log(2/\delta_j)}}{n_j\eps_jh^d}\xi^{(j)}(\cdot)
		\text{ for }
		j=0,1,\dots,m.
	\end{equation}
	We then release
	\begin{equation}\label{eq:pvt-Tj}
		\left\{T^{(j)}_h(x_0)+\tilde{\xi}^{(j)}_h(x_0)\right\}_{j=0}^m.
	\end{equation}

	\noindent The next proposition asserts that the transcripts described above matches the required distributed privacy requirements, and its proof follows by results from \cite{hall2013differential}.
	
	\begin{proposition}\label{pr:privacy-guarantee}
		For any $h\in [0,1]$ the transcripts $\{T^{(j)}_h(x_0)+\tilde{\xi}^{(j)}_h(x_0):0\le j\le m\}$ described above satisfies $(\eps_j,\delta_j)$ differential privacy distributed across servers $j\in\{0,1\dots,m\}$.
	\end{proposition}
	
	The optimal bandwidth choice $h_{opt}$ is given by the solution to \eqref{eq:optimal-bandwidth-choice}. To account for the additional $\delta$ factor for approximate privacy, we now define $h_{opt,\delta}$ which is the solution to:
	\begin{equation}\label{eq:optimal-h-delta}
		(n_0\wedge n_0^2\eps_{0}^2r^d)r^{2\beta+d}
		+
		\displaystyle\sum_{j=1}^m
		(n_{j}\wedge n_{j}^2\eps_{j}^2r^d)
		r^{2\beta\gamma_j+d}
		=~\log\left(
		\dfrac{2}{\delta_{\min}}
		\right)
	\end{equation}
	Let us now define the weights
	\begin{align*}\label{eq:upbd-wts}
		v_{j}=(n_j\wedge n_j^2\eps_j^2h_{opt,\delta}^d)h_{opt,\delta}^{\gamma_j\beta}
		\quad 
		&\text{for}
		\quad
		j=0,\dots,m
		\\
		u_{j}=\dfrac{v_{j}}{\sum_{j=0}^mv_{P}}
		\quad 
		&\text{for}
		\quad
		j=0,\dots,m
		.\numberthis
	\end{align*}
	With these weights, the target server computes a weighted average of $\{T^{(j)}_h(x_0)+\tilde{\xi}_h^{(j)}(x_0):0\le j\le m\}$ as follows:
	\begin{align*}\label{eq:def-Tn}
		\tilde{T}_h(x_0) := 
		&~
		u_{0}
		\left(
		T^{(0)}_h(x_0)
		+
		\tilde{\xi}^{(0)}(x_0)
		\right) 
		+
		\sum_{j=1}^m u_{j}
		\left(
		T^{(j)}_h(x_0)
		+\tilde{\xi}^{(j)}(x_0)
		\right)
		.\numberthis
	\end{align*}
	Finally our classifier is given by
	\begin{equation}\label{eq:def-fhat}
		\hat{f}(x_0):=
		\bbm{1}(\tilde{T}_{h_{opt,\delta}}(x_0)\ge 0)
	\end{equation}
	where $h_{opt,\delta}$ is the solution to~\eqref{eq:optimal-h-delta}. The following theorem provides an upper bound for the excess risk of this classifier.
	\begin{theorem}\label{th:up-bd}
		Let $r$ be the solution to \eqref{eq:optimal-bandwidth-choice}. Let $\hat{f}$ be the classifier defined in~\eqref{eq:def-fhat} based on the weighted kernel estimator~\eqref{eq:def-Tn}. Then,
		$$
		\sup_{(P_0,\ldots,P_m) \in \Pi(\alpha,\beta,\bm \gamma,\mu)}
		\calE_0(\hat{f})
		\le C_*
		r^{\beta(1+\alpha)}
		\left(
		\log\left(
		\dfrac{1}{\delta_{\min}}
		\right)
		\right)^{\frac{\beta(1+\alpha)}{2\beta\gamma_{\min}+d}}
		$$
		where $C_*$ is a constant depending on $L,d,\alpha,\beta,\gamma_j$, while   $\gamma_{\min}=\min\{1,\gamma_1,\dots,\gamma_m\}$ and $\delta_{\min}=\min\{\delta_0,\dots,\delta_m\}$.
	\end{theorem}
	
	%%%%%%%%%%%%%%%%%%%%%%%%%%%%%%%%%%%%%%%%%%
	\subsection{Minimax Lower Bounds}
	%%%%%%%%%%%%%%%%%%%%%%%%%%%%%%%%%%%%%%%%%%
	
	The above theorem bounds the error rate of our kernel based classifier. Alongside the upper bound above, in this subsection we derive the minimax lower bound on the excess risk, to establish that our kernel based classifier is minimax optimal up to logarithmic factors.

	We introduce a general data processing inequality which extends the findings presented in \cite{Cai2023FL-NP-Regression}. This new result provides a bound on the total variation (TV) distance between the push forward measures of the transcripts $\P_\sigma^{\tilde{T}}$ and $\P_{\sigma'}^{\tilde{T}}$, utilizing the TV distance between their underlying distributions. Such an inequality may be of independent interest beyond the current setting due to its broader applicability.
	
	\begin{lemma}\label{lemma:total-variation-two-point}
		For any subset $\calS \subseteq \{0,\ldots,m\}$, the TV distance is bounded as follows:
		\begin{equation}\label{eq:total_variation_split}
			\mathrm{TV}\left( \PP_{\sigma}^{\tilde{T}}, \PP_{\sigma'}^{\tilde{T}} \right) \leq \sqrt{2} \sqrt{\sum_{j \in S} \bar{\varepsilon}_j \left( e^{\bar{\varepsilon}_j} - 1 \right) + \sum_{j \in S^c} n_j \mathrm{KL}( P_{j,\sigma}, P_{j,\sigma'})} + 4 \sum_{j \in S} e^{\bar{\varepsilon}_j} n_j \delta_j \rho_j,
		\end{equation}
		where $\bar{\varepsilon}_j = 6n_j\varepsilon_j \rho_j$ and $\rho_j = \mathrm{TV}(P_{j,\sigma}, P_{j,\sigma'})$.
	\end{lemma}
	
	A notable aspect of this lemma, in contrast to previous results in \cite{Cai2023FL-NP-Regression}, is the inclusion of interactions between the target and source servers. Consequently, all the information about the transcripts from the source servers is encapsulated within the final transcript ${\tilde{T}}$ computed by the target server. This approach allows for heterogeneity not only in sample sizes and privacy budgets but also in the data distributions across different servers.
	
	In constructing the lower bound, we formulate a family of distributions $P_{j,\sigma}$ for each $j = 1, \ldots, m$, where $\sigma$ is a vector of $\{+1, -1\}^M$ and $M$ is suitably chosen. We apply Assouad's Lemma to bound the total variation distance between the push forward measures of the transcripts $\P_\sigma^{\tilde{T}}$ and $\P_{\sigma'}^{\tilde{T}}$ for $\sigma$ and $\sigma'$ differing in one entry, as stipulated by Lemma \ref{lemma:total-variation-two-point}. The rest of the construction of the lower bound crucially depends on the selection of the set $\calS$. This choice corresponds to the index of servers where the privacy cost significantly outweighs the non-private risk. The major contributions to the upper bound for the TV distance are $\bar{\varepsilon}_j \left( e^{\bar{\varepsilon}_j} - 1 \right)$ from the privacy-stringent servers, and the KL divergence for the samples on non-stringent servers, $n_j \mathrm{KL}(P_{j,\sigma}, P_{j,\sigma'})$. The optimal balance of contributions is achieved when $\calS = \{j : \varepsilon_j \leq (r^dn_j)^{-1/2}\}$.
	
	The following theorem establishes the fundamental cost of privacy for the nonparametric classification problem in the distributed privacy setting.
	\begin{theorem}\label{th:lo-bd}
		Suppose $\delta_j$'s are such that $\sum_j  n_j\delta_j = o(1)$, then there exists a $c >0$ not depending on $n_j$ for $j =0,\ldots,m$ such that
		$$
		\inf_{\widehat{f}\in \cM(\bm\varepsilon,\bm\delta)}\sup_{(P_0,\ldots,P_m) \in \Pi} 
		\cE_0(\widehat f) \geq cr^{\beta(1+\alpha)}
		$$
		where $r$ is the solution to~\eqref{eq:optimal-bandwidth-choice}.
	\end{theorem}
	A comparison with the upper bound on the excess risk obtained from Theorem~\ref{th:up-bd} now establishes the rate optimality, up to a logarithmic factor, of our kernel-based differentially private distributed classifier.
	
	%%%%%%%%%%%%%%%%%%%%%%%%%%%%%%%%%%%%%%%%%%
	\section{Data-driven Adaptive Classifier}\label{sec:adaptation}
	%%%%%%%%%%%%%%%%%%%%%%%%%%%%%%%%%%%%%%%%%%
	
	In practice, the smoothness and transfer exponent parameters are unknown, making it challenging to select the correct bandwidth $h$. To address this, we will use an estimator based on the Lepski method to choose $h$ from a grid of possible values. While choosing the exact optimal $h$ is infeasible, this method adapts to the unknown parameters $\beta$ and $\gamma_1,\dots,\gamma_m$. 
	
	This section is divided into two subsections. First, we consider the transfer homogeneous case, where the source populations have the same relative transfer exponent $\gamma_j=\gamma$ for $j=1,\dots,m$  with respect to the target. Here we allow $m$ to grow at an appropriately slow rate as $n$ increases. In the second subsection, we address the general, heterogeneous case with transfer exponents $(\gamma_1,\dots,\gamma_m)$ with the important restriction that the number of sources $m$ is finite.
	
	To choose the best candidate bandwidth, we define a grid of possible choices for $h$ as:
	$$
	\calH=\{2^{-j}:j=0,1,\dots,(\log n_*)/d\}, 
	\text{ where }
	n_* = \sum_{j=0}^m n_j\wedge n_j^2\eps_j^2.
	$$
	Let $\Delta^m=\{w:w_i\in [0,1],\sum_{i=0}^mw_i=1\}$ denote the $m$-dimensional simplex. For a weight vector $w=(w_0,w_1,\dots,w_m)\in\Delta^m$  we define 
	\begin{equation}\label{eq:def-Tx0-h-w}
		\tilde{T}(x_0,h,w)
		:=
		\sum_{j=0}^mw_jT^{(j)}_h(x_0)
		+
		\sum_{j=0}^mw_j
		\dfrac{\sqrt{2c_K\log(2|\calH|/\delta_j)}|\calH|}{n_j\eps_jh^d}
		{\xi}^{(j)}(x_0)
		\text{ for }
		h,w\in[0,1]
	\end{equation}
	where $T_h^{(j)}(\cdot)$ is as defined in \eqref{eq:def-Tnj-np} and ${\xi}^{(j)}(\cdot)$ are independent mean zero Gaussian processes with covariance kernel $K(\cdot/h)$. 
	
	%%%%%%%%%%%%%%%%%%%%%%%%%%%%%%%%%%%%%%%%%%
	\subsection{Adaptation under Source Homogeneity}\label{sec:adapt-homog}
	%%%%%%%%%%%%%%%%%%%%%%%%%%%%%%%%%%%%%%%%%%
	
	In this subsection we consider the sources to have transfer homogeneity, i.e., every source has identical transfer exponent $\gamma_j=\gamma$ for $j=1,\dots,m$. It is then intuitive to weigh the estimators $T_h^{(j)}$ and the noise $\xi^{(j)}$ with the weights proportional to $n_j\wedge n_j^2\eps_j^2h^d$ for the sources. More specifically, this gives the restricted set of weights:
	\begin{equation}\label{eq:def-calW}
		\calW(h)
		:=
		\bigg\{
		(w_0,w_1,\dots,w_m)
		:~
		w_0\in[0,1],\,		
		w_j=\dfrac{(1-w_0)u_j(h)}{\sum_{j=1}^mu_j(h)}
		\text{ for }
		j=1,\dots,m
		\bigg\}%\numberthis
	\end{equation}
	where $u_j(h):=n_j\wedge n_j^2(\eps_j/|\calH|)^2h^d$. Note that the privacy requirements dictate that $\calW$ depends on $h$. The rationale behind this set of weights comes from \eqref{eq:upbd-wts} restricted to the transfer homogeneous setting. We determine the deviation of $\tilde{T}(x_0,h,w)$ around its expectation through
	\begin{align*}\label{eq:var-ada-homog}
		v_0(h,w)
		=&~
		w_0^2\left(
		\dfrac{c_Kg_{\max}}{3n_0h^d}+
		\dfrac{2c_K^2\log(2|\calH|/\delta_0)}{n_0^2(\eps_0/|\calH|)^2h^{2d}}\right)\\
		&~
		+
		(1-w_0)^2
		\left(
		\sum_{j=1}^m
		\dfrac{(u_j(h))^2}{
			(\sum_{j=1}^mu_j(h))^2}
		\bigg(
		\dfrac{c_Kg_{\max}}{3n_jh^d}
		+
		\dfrac{2c_K^2\log(2|\calH|/\delta_j)}{n_j^2(\eps_j/|\calH|)^2h^{2d}}
		\bigg)\right)
		\numberthis.
	\end{align*} 
	The above metric measures the scale of noise present in $\tilde{T}(x_0,h,w)$. We therefore define the signal to noise ratio as
	\begin{equation}\label{eq:def-rho0-h}
		\hat{\rho}_0(h)=\max_{w\in\calW(h)}\dfrac{(\tilde{T}(x_0,h,w))^2}{v_0(h,w)}.
	\end{equation}
	Then the adaptive choice of $h$ in the transfer homogeneous setting is given by
	$$
	h_{0}
	=\begin{cases}
		\min\left\{h\in\calH:\hat{\rho}_0(h)>
		4.5
		\log(2n_*|\calH|)
		\right\}
		\quad &\text{if }
		\underset{h\in\calH}{\max}\,\hat{\rho}_0(h)
		>
		4.5
		\log(2n_*|\calH|)\\
		{\rm argmax}_{h}\hat{\rho}_0(h)\quad&\text{otherwise}.
	\end{cases}
	$$
	Define 
	\[
	w_{*(0)}={\rm argmax}_{w\in\calW(h_0)}\dfrac{(\tilde{T}(x_0,h_0,w))^2}{v_0(h_0,w)}.
	\]
	The adaptive classifier is now defined as
	\begin{equation}\label{eq:def-f-ada0}
		\hat{f}_{0}(x):=\bbm{1}(\tilde{T}(x,h_0,w_{*(0)})>0).
	\end{equation}
	Algorithm~\ref{alg:ada-homog} summarizes the steps of the adaptive procedure.

	\begin{algorithm}[h]
		\caption{Data-adaptive mechanism for bandwidth selection under Source Homogeneity}\label{alg:ada-homog}
		\begin{algorithmic}[1]
			\State \textbf{Input:} test point: $x_0\in[0,1]^d$; data: $\{(X_i^{(j)},Y_i^{(j)}):1\le i\le n_j, j\in\{0,1,\dots,m\}\}$; privacy parameters: $\{(\eps_j,\delta_j):j\in\{0,1,\dots,m\}\}$; kernel: $K(\cdot)$.
			\State Compute the grid of bandwidths: $\calH=\{2^{-j}:j=0,1,\dots,\floor{\log(n_*)/d}\}$.
			\For{$h$ in $\calH$}
			\State Compute the set of weights $\calW(h)$ following \eqref{eq:def-calW}.	
			\State Compute the signal-to-noise ratio:
			$$
			\hat{\rho}_0(h)=\max_{w\in\calW(h)}\dfrac{(\tilde{T}(x_0,h,w))^2}{v_0(h,w)}
			$$
			where $\tilde{T}(x_0,h,w)$ and $v_0(h,w)$ are as defined in \eqref{eq:def-Tx0-h-w} and \eqref{eq:var-ada-homog} respectively.
			\EndFor
			\State Choose the bandwidth as
			$$h_{0}
			=\begin{cases}
				\min\left\{h\in\calH:\hat{\rho}_0(h)>
				4.5
				\log(2n_*|\calH|)
				\right\}
				\quad &\text{if }
				\underset{h\in\calH}{\max}\,\hat{\rho}_0(h)
				>
				4.5
				\log(2n_*|\calH|)\\
				{\rm argmax}_{h\in\calH}\hat{\rho}_0(h)\quad&\text{otherwise}.
			\end{cases}
			$$
			\State Choose the best weight as
			$$
			w_{*(0)}={\rm argmax}_{w\in\calW(h_0)}\dfrac{(\tilde{T}(x_0,h_0,w))^2}{v_0(h_0,w)}.
			$$
			\State \textbf{Output:} $\hat{f}_{0}(x):=\bbm{1}(\tilde{T}(x,h_0,w_0)>0)$.
		\end{algorithmic}
	\end{algorithm}
	The following theorem states the excess risk of the adaptive classifier in terms of the regression function parameter $\alpha,\beta$, the transfer exponent $\gamma$ and the privacy constraints.
	
	\begin{theorem}\label{th:adapt-homog}
		Let $r$ be the solution to \eqref{eq:optimal-bandwidth-choice} with $\gamma_j=\gamma$ for $j=1,\dots,m$. Let $\hat{f}_0$ be the data adaptive classifier defined in~\eqref{eq:def-f-ada0}. Then,
		$$
		\sup_{(P_0,\ldots,P_m) \in \Pi(\alpha,\beta,\gamma,\mu)}
		\calE_0(\hat{f}_0)
		\le 
		C_*'
		r^{\beta(1+\alpha)}
		\left[(
		\log(n_*|\calH|)
		\log\left(2|\calH|/\delta_{\min}
		\right))^{\frac{\beta(1+\alpha)}{2\beta(1\wedge\gamma)+d}}\vee |\calH|^{\frac{2\beta(1+\alpha)}{d}}
		\right]
		$$
		where $C_*'$ is a constant depending on $m,L,d,\alpha,\beta,\gamma$, while $\delta_{\min}=\min\{\delta_0,\dots,\delta_m\}$.
	\end{theorem}
	
	It is instructive to compare the above rate with the rate from Theorem~\ref{th:up-bd}. Since $|\calH|={\rm polylog}(n_*)$ it follows that the excess risk of the adaptive estimator is worse by a multiplicative factor of ${\rm polylog}(n_*))$. This is due to two reasons. Firstly, a factor of $\log n_*$ is expected as a cost of adaptation, and can be found in the transfer learning classification setup considered in \cite{cai2021transfer}. Secondly, to ensure that the $(\eps_j,\delta_j)$ privacy requirements are satisfied throughout the adaptation procedure, we require that for every $h\in\calH$ our estimators are $(\eps_j/|\calH|,\delta_j/|\calH|)$ differentially private. This contributes an extra factor of $|\calH|$ to the rate obtained in Theorem~\ref{th:adapt-homog}. 
	
	An important special case of the above theorem is the server homogeneous case, where sample sizes and the privacy parameters are the same for every server, i.e., $n_j=n$, $\eps_j=\eps$ and $\delta_j=\delta$ for $j=0,1,\ldots,m$. The following corollary  describes this special case.
	
	\begin{corollary}\label{cor:adapt-homog}
		Let $r$ be the solution to \eqref{eq:optimal-bandwidth-choice} with $n_j =n,\varepsilon_j =\varepsilon,\delta_j = \delta$ and $\gamma_j =\gamma$ for all $j=1,\ldots,m$. Let $\hat{f}_0$ be the data adaptive classifier defined in~\eqref{eq:def-f-ada0}. Then,
		\begin{align*}
			\sup_{(P_0,\ldots,P_m) \in \Pi(\alpha,\beta,\gamma,\mu)}
			\cE_0(\hat{f}_0)
			\le 
			C_*'\bigg[ L_N^{(ada)}\bigg\{
			&~
			\left( n_0^{\frac{1}{2\beta + d}} \wedge (n_0^2\varepsilon^2_0)^{\frac{1}{2\beta + 2d}}\right) \\
			&~+ \left( (mn)^{\frac{1}{2\beta \gamma+ d}} 
			\wedge (mn^2\varepsilon^2)^{\frac{1}{2\beta \gamma+ 2d}}\right) \bigg\}^{-\beta(1+\alpha)}
			\wedge 1\bigg]
		\end{align*}
		where $L_N^{(ada)}$ is given by $\left[(
		\log(n_*|\calH|)
		\log\left(2|\calH|/\delta
		\right))^{\frac{\beta(1+\alpha)}{2\beta(1\wedge\gamma)+d}}\vee |\calH|^{\frac{2\beta(1+\alpha)}{d}}
		\right]$, and $C_*'$ is a constant depending on $m,L,d,\alpha,\beta,\gamma$.
	\end{corollary}

	A comparison with Theorem~\ref{th:homogeneous-minimax-rate} shows that this rate is minimax optimal up to a factor polynomial in logarithmic terms.
	
	%%%%%%%%%%%%%%%%%%%%%%%%%%%%%%%%%%%%%%%%%%
	\subsection{General Adaptation for Multiple Sources}\label{sec:adapt-gen}
	%%%%%%%%%%%%%%%%%%%%%%%%%%%%%%%%%%%%%%%%%%
	
	We now shift to the general setting where we no longer constrain $\gamma_1,\dots,\gamma_m$ to be all equal. Note that for optimal estimation (as in Theorem~\ref{th:up-bd}), one requires knowledge of potentially $m$ many different parameters $\gamma_1,\dots,\gamma_m$. The adaptation procedure therefore requires optimizing over all possible weights in $w\in\Delta^m$. When $m$ increases with $n$, the adaptation to this growing number of parameters necessarily worsens the rate of decay for the excess risk. We will not delve further into issue and focus instead on the case where $m$ is finite and does not increase with $n$. 
	
	In this general case we must consider all possible weight vectors $w\in\Delta^m$. As in \eqref{eq:var-ada-homog} earlier we compute an approximate variance of $\tilde{T}(x_0,h,w)$ as 
	\begin{equation}\label{eq:var-ada-gen}
		v(h,w) = 
		\sum_{j=0}^m
		w_j^2
		\left(
		\dfrac{c_Kg_{\max}}{3n_jh^d}
		+
		\dfrac{2c_K^2\log(2|\calH|/\delta_j)}{n_j^2(\eps_j/|\calH|)^2h^{2d}}
		\right)
		.
	\end{equation}
	Let us define the signal-to-noise ratio index $\hat{\rho}(h)$:
	$$
	\hat{\rho}(h)=\max_{w\in\Delta^m}\dfrac{(\tilde{T}(x_0,h,w))^2}{v(h,w)},
	$$
	In this general setting, the adaptive choice of $h$ is given by
	$$
	h_*
	=\begin{cases}
		\min\left\{h\in\calH:\hat{r}(h)>
		C_*
		\log(n_*|\calH|(m+1))
		\right\}
		\quad &\text{if }
		\underset{h\in\calH}{\max}\,\hat{\rho}(h)
		>
		C_*
		\log(n_*|\calH|(m+1))\\
		{\rm argmax}_{h}~\hat{\rho}(h)\quad&\text{otherwise}.
	\end{cases}
	$$
	where $C_*=2.25(m+1)$.
	Defining
	$$
	w_*={\rm argmax}_{w\in\Delta^m}\dfrac{(T(x_0,h_*,w))^2}{v(h_*,w)}.
	$$
	we obtain the adaptive classifier
	\begin{equation}\label{eq:def-f-ada}
		\hat{f}_a(x):=\bbm{1}(T(x,h_*,w_*)>0).
	\end{equation}
	Similar to Algorithm~\ref{alg:ada-homog} we have Algorithm~\ref{alg:ada-gen} which presents the adaptive procedure allowing for possible heterogeneity among servers.
	
	\begin{algorithm}[h]
		\caption{Data-adaptive mechanism for bandwidth selection in the general case}\label{alg:ada-gen}
		\begin{algorithmic}[1]
			\State \textbf{Input:} test point: $x_0\in[0,1]^d$; data: $\{(X_i^{(j)},Y_i^{(j)}):1\le i\le n_j, j\in\{0,1,\dots,m\}\}$; privacy parameters: $\{(\eps_j,\delta_j):j\in\{0,1,\dots,m\}\}$; kernel: $K(\cdot)$.
			\State Compute the grid of bandwidths: $\calH=\{2^{-j}:j=0,1,\dots,\floor{\log(n_*)/d}\}$.
			\For{$h$ in $\calH$}
			\State Compute the signal-to-noise ratio:
			$$
			\hat{\rho}(h)=\max_{w\in \Delta^m}\dfrac{(\tilde{T}(x_0,h,w))^2}{v(h,w)}
			$$
			where $\tilde{T}(x_0,h,w)$ and $v(h,w)$ are as defined in \eqref{eq:def-Tx0-h-w} and \eqref{eq:var-ada-gen} respectively.
			\EndFor
			\State For $C_*=2.25(m+1)$ choose the bandwidth as
			$$h_*
			=\begin{cases}
				\min\left\{h\in\calH:\hat{\rho}(h)>
				C_*
				\log(2n_*|\calH|)
				\right\}
				\quad &\text{if }
				\underset{h\in\calH}{\max}\,\hat{\rho}(h)
				>
				C_*
				\log(2n_*|\calH|)\\
				{\rm argmax}_{h\in\calH}\hat{\rho}(h)\quad&\text{otherwise}.
			\end{cases}
			$$
			\State Choose the best weight as
			$$
			w_*={\rm argmax}_{w\in \Delta^m}\dfrac{(\tilde{T}(x_0,h_0,w))^2}{v(h_0,w)}.
			$$
			\State \textbf{Output:} $\hat{f}_{a}(x):=\bbm{1}(\tilde{T}(x,h_*,w_*)>0)$.
		\end{algorithmic}
	\end{algorithm}
	
	The following theorem verifies the efficacy of the general adaptive procedure.
	
	\begin{theorem}\label{th:adapt}
		Let $r$ be the solution to \eqref{eq:optimal-bandwidth-choice}. Let $\hat{f}_a$ be the data adaptive classifier defined in~\eqref{eq:def-f-ada}. Then,
		$$
		\sup_{(P_0,\ldots,P_m) \in \Pi(\alpha,\beta,\bm \gamma,\mu)}
		\calE_0(\hat{f})
		\le 
		C_*'
		r^{\beta(1+\alpha)}
		\left[(
		\log(n_*|\calH|)
		\log\left(2|\calH|/\delta_{\min}
		\right))^{\frac{\beta(1+\alpha)}{2\beta\gamma_{\min}+d}}\vee |\calH|^{\frac{2\beta(1+\alpha)}{d}}
		\right]
		$$
		where $C_*'$ is a constant depending on $m,L,d,\alpha,\beta,\gamma_j$, while $\gamma_{\min}=\min\{1,\gamma_1,\dots,\gamma_m\}$ and $\delta_{\min}=\min\{\delta_0,\dots,\delta_m\}$.
	\end{theorem}
	
	\begin{remark}{\rm
			Oftentimes users might find it helpful to restrict the set of weights based on prior knowledge. For example, under source homogeneity we allowed the restricted set $\calW(h)$. Similarly it is possible to use ad-hoc choices of weights $w$, based for example, on the sample proportions for each server. Our algorithm~\ref{alg:ada-gen} automatically allows these specific choices of weights and can be used to select the bandwidth $h$ to be used with the classifiers weighted with respect to $w$.
	}\end{remark}

	%%%%%%%%%%%%%%%%%%%%%%%%%%%%%%%%%%%%%%%%%%
	\section{Numerical Studies}\label{sec:numerical-studies}
	%%%%%%%%%%%%%%%%%%%%%%%%%%%%%%%%%%%%%%%%%%
	
	The data-driven classifier proposed in this paper is easy to implement. In this section, we examine the numerical performance of our methods through various simulated and real data experiments. In the first subsection, we conduct simulation studies to compare our proposed classifier with alternative methods across different parameter settings and varying levels of problem difficulty. In the second subsection, we evaluate the performance of the proposed classifier on a real dataset, verifying its practical effectiveness.

	%%%%%%%%%%%%%%%%%%%%%%%%%%%%%%%%%%%%%%%%%%
	\subsection{Simulation Study}
	%%%%%%%%%%%%%%%%%%%%%%%%%%%%%%%%%%%%%%%%%%
	
	We compare the prediction accuracy of our classifier with other existing methods in the literature across various parameter settings. 
	
	%%%%%%%%%%%%%%%%%%%%%%%%%%%%%%%%%%%%%%%%%%
	\subsubsection{Simulation Design}\label{sec:simulation-design}
	%%%%%%%%%%%%%%%%%%%%%%%%%%%%%%%%%%%%%%%%%%
	
	In our simulation study, we consider a setup with $m$ source servers and a single target server. All source servers are assumed to have the same data quality. Additionally, for most of the simulation setups we assume that both the source and target servers have an equal number of observations, denoted by $n$. The privacy budget for each server is given by $(\varepsilon, \delta)$, where we set $\delta = n^{-2}$.
	
	The data for both target and source servers is generated as follows: the marginal distribution for $\bm{X}$ on both the target and source servers is supported on a two-dimensional cube, \([0, 1]^2\), with a uniform distribution over its support.
	
	The conditional distribution for the target server, given $\bm{X}$, is defined as:
	\[
	P_T(Y = 1 | \bm{X}) = \eta_T(\bm{X}) = 
	1\wedge   
	\left( 
	\frac{1}{2} + \mathrm{sign}\left(
	\bigg(x_1 - \frac{1}{2}\bigg)\bigg(x_2 - \frac{1}{2}\bigg)\right) 
	\left|x_1 - \frac{1}{2}\right|^{\frac{1}{4}} \left|x_2 - \frac{1}{2}\right|^{\frac{1}{4}}
	\right)_+
	\]
	where for $a\in\RR$ we write $a_+=\max\{a,0\}$.	For the source servers, the conditional distribution is defined as:
	\[
	P_S(Y = 1 | \bm{X}) = 
	1\wedge 
	\left( 
	\frac{1}{2} + \mathrm{sign}\left(\eta_T(\bm{X}) - \frac{1}{2}\right) \left|\eta_T(\bm{X}) - \frac{1}{2}\right|^{\gamma}\right)_+.
	\]
	
\begin{figure}[h]
		\centering
		\includegraphics[scale = 0.6]{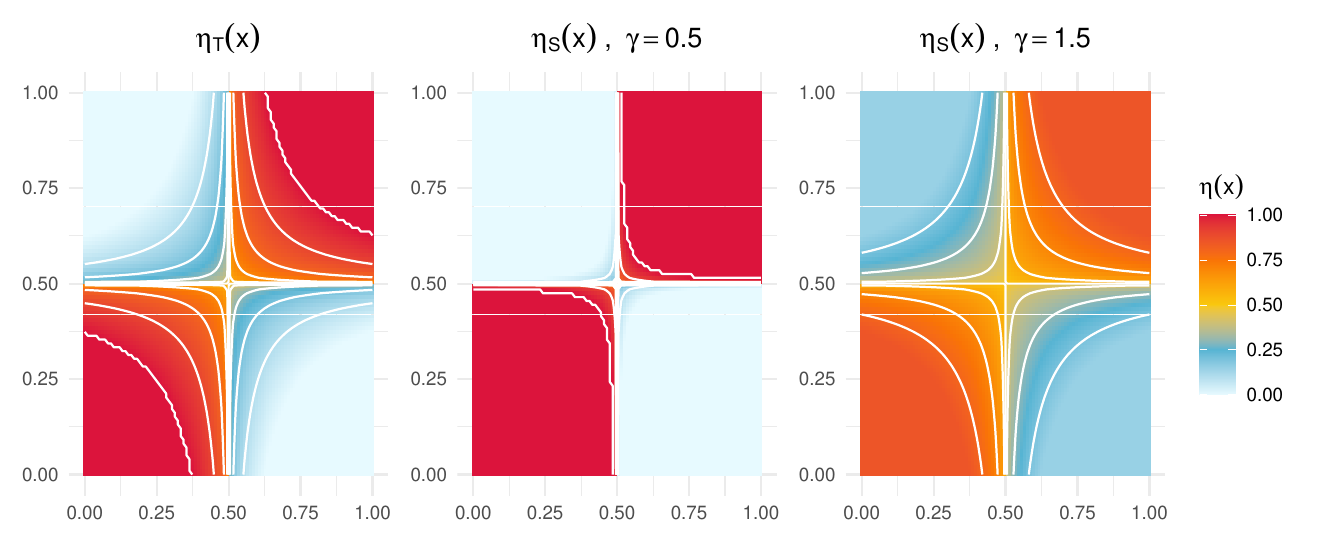}
		\caption{Regression functions \(\eta_T\) and \(\eta_S\) for \(\gamma \in \{0.5, 1.5\}\).}
		\label{fig:reg_func}
	\end{figure}

	The proposed construction for these regression functions satisfies key assumptions, such as smoothness and margin conditions. For a visual representation of these regression functions, see Figure \ref{fig:reg_func}.
	
	We would be comparing the following methods:
	\begin{itemize}
		\item \textbf{Distributed Transfer Learning with kernels (DTK)} is the proposed method with two different choices of kernel, the gaussian kernel and the triangular kernel. We select the best tuning parameters $h$ the bandwidth and, $w \in[0,1]$ the relative weight of the target data w.r.t the source data. The parameter grid for those are $w \in \left\{\frac{i}{100} \,:\, i\in\{0,\ldots,100\}\right\}$ and $h \in \{2^{-i} \,:\, i\in\{1,\ldots,7\}\}$.
		\item \textbf{Distributed Transfer Learning using histogram (DT-HIST)} is a variant of the method proposed in \cite{berrett2019classification} adapted to the transfer learning setting with distributed privacy. This method also has the the two parameters $w$ and $h$, where $h$ denotes the histogram bin size for this estimator. We vary both the parameters on the same grid as \textbf{DTK}.
		\item \textbf{Adaptive Distributed Learning with kernels (AdaptDTK)} This is the adaptive version of the DTK classifier, proposed in Section \ref{sec:adaptation}, which automatically tunes the hyper-parameters $h$ and $w$ based on a Lepski-style method. 
	\end{itemize}
	
	%%%%%%%%%%%%%%%%%%%%%%%%%%%%%%%%%%%%%%%%%%
	\subsubsection{Effect of Source Data}\label{sec:effect_of_source}
	%%%%%%%%%%%%%%%%%%%%%%%%%%%%%%%%%%%%%%%%%%
	
	In this section, we examine the impact of source data on classification performance. We compare our transfer learning approach to a method that relies solely on target data. The latter method, referred to as \textbf{targetDTK}, serves as our baseline for comparison.
	
	To assess the utility of source data, we plot the best accuracy obtained from each method across a range of hyperparameters, as detailed in Section \ref{sec:simulation-design}. This analysis provides insight into the benefits of incorporating source data in our classification process.
	\begin{figure}[h]
		\centering
		\includegraphics[scale = 0.67]{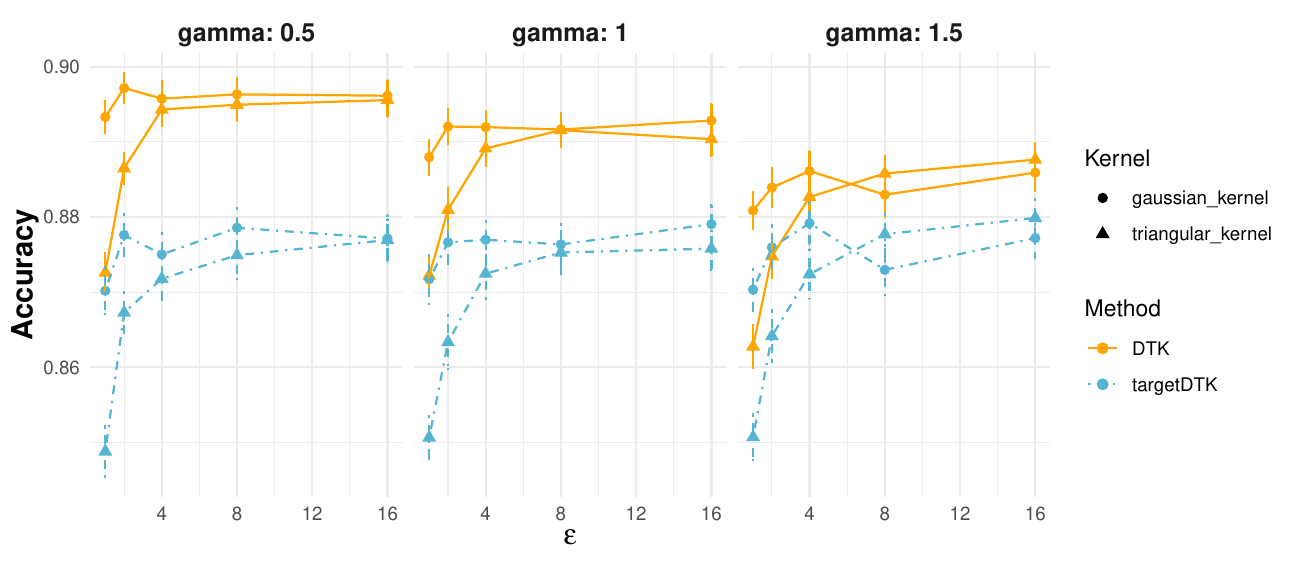}
		\caption{Effect of source data against $\varepsilon$ for \(\gamma \in \{0.5, 1, 1.5\}\)}
		\label{fig:effect_of_source}
	\end{figure}
	In Figure \ref{fig:effect_of_source}, we contrast our transfer learning method with a naive approach that doesn't leverage source data. To ensure fairness, we limit the scenario to just one additional source server ($m = 1$). Introducing multiple source servers would widen the performance gap between transfer learning methods and non-transfer learning methods. Also we set the number of observations on the source and target server to be $n =100$, the reason for choosing a comparatively smaller $n$ is for better visualization purposes, for larger $n$ we expect the same phenomena , although the gain in accuracy would be difficult to perceive visually.
	
	Within this setup, we vary two problem parameters, $\varepsilon$ and $\gamma$. Across various kernel choices and values of $\varepsilon$ and $\gamma$, our transfer learning method consistently outperforms the alternatives. As privacy constraints relax (i.e., $\varepsilon$ increases), performance improves across all methods. Notably, as $\gamma$ increases, the performance advantage gained from using source servers diminishes. This phenomenon is expected: with higher $\gamma$, the quality of source data decreases, consequently reducing the performance gain. Since it is always beneficial to use source data we do not compare the target only method in the further sections.
	
	%%%%%%%%%%%%%%%%%%%%%%%%%%%%%%%%%%%%%%%%%%
	\subsubsection{Comparison of Our Adaptive Classifiers with Other Methods}
	%%%%%%%%%%%%%%%%%%%%%%%%%%%%%%%%%%%%%%%%%%
	
	In this section we compare our adaptive classifier \textbf{AdaptDTK} to non-adaptive methods like \textbf{DTK} and \textbf{DT-HIST} across different $\varepsilon$, $\gamma$ and $m$ values. It is important to note that both non-adaptive methods are tuned with optimal hyperparameters based on test data, resulting in a potential performance advantage over the adaptive method. 
	
	However, our results show that the ``cost of adaptation"---the difference in accuracy between the adaptive method and the non-adaptive methods---is relatively small. Despite the adaptive method's lack of knowledge about oracle hyper-parameters, its performance is only marginally lower compared to methods with oracle tuning. This indicates that our adaptive approach is effective, even without precise hyperparameter tuning.
	
	\paragraph{Effect of Privacy Budget:} 
	
	\begin{figure}[h]
		\centering
		\includegraphics[scale = 0.7]{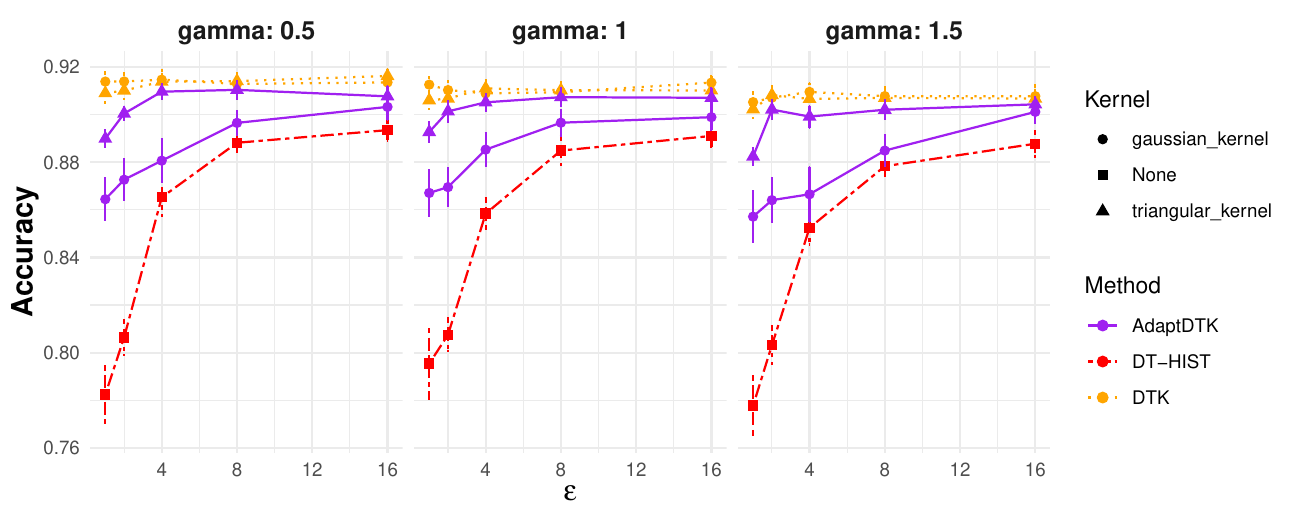}
		\caption{Accuracy v/s $\varepsilon$ for \(\gamma \in \{0.5, 1, 1.5\}\)}
		\label{fig:effect_of_privacy}
	\end{figure}
	Here, we maintain the same experimental setup as described in Section \ref{sec:effect_of_source} but with larger $n =500$. Our kernel-based method, \textbf{AdaptDTK}, along with \textbf{DTK} (with oracle hyper-parameter tuning), consistently outperforms \textbf{DT-HIST}. Particularly in scenarios with low privacy budgets, the performance gap is significant, indicating that privatizing kernel-based methods empirically yields better results compared to privatizing histogram bin count-based methods.
	
	Regarding the choice of kernel, we note that the triangular kernel demonstrates exceptional adaptability, showing only minimal performance degradation compared to its oracle-tuned counterpart. This difference in performance is likely attributed to the fact that our theoretical adaptive procedure was developed with bounded kernels (like triangular kernels) in mind. 
	
	\paragraph{Effect of distributed data:} 
	\begin{figure}[h]
		\centering
		\includegraphics[scale = 0.7]{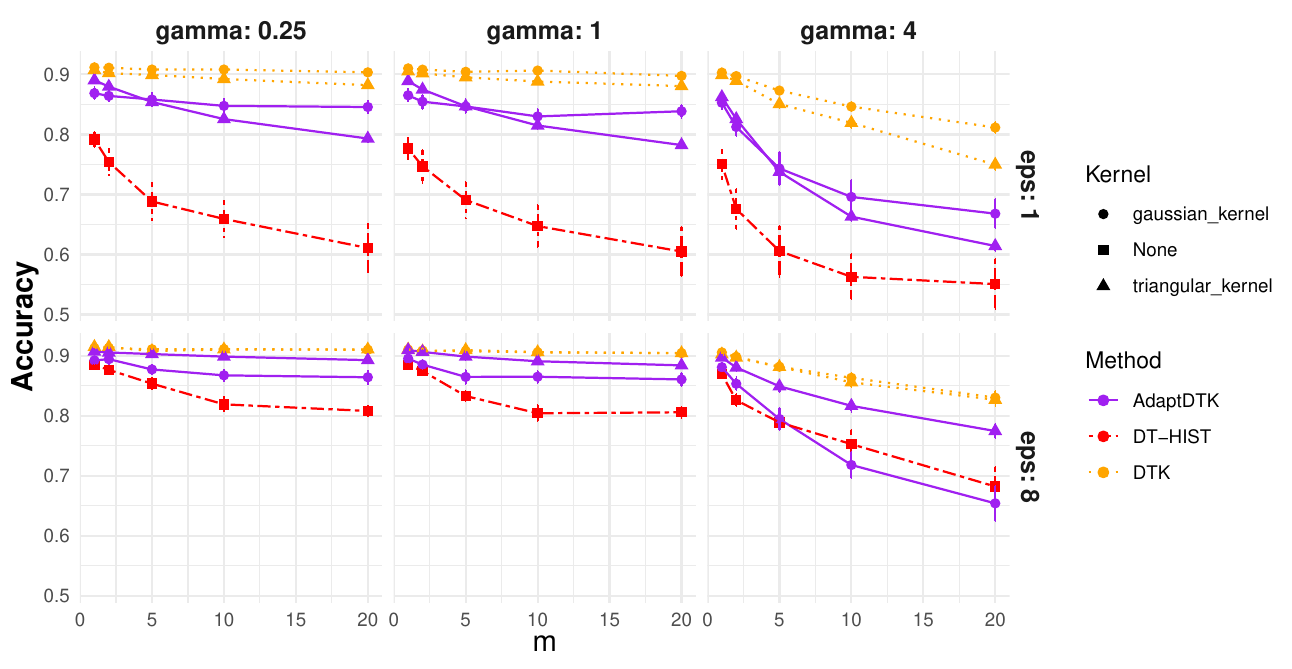}
		\caption{Accuracy v/s $m$ for \(\gamma \in \{0.25, 1, 4\}\)}
		\label{fig:effect_of_m}
	\end{figure}
	In our experimental setup, we maintain a constant total sample size of $500$ across both target and source servers. By varying the number of source servers $m$ from $1$ to $20$, we delve into the impact of data distribution under a fixed total sample size. Additionally, we explore how the source data quality, represented by $\gamma$, and the privacy budget, $\varepsilon$, influence the outcomes. Across all methods, a consistent pattern emerges: as data becomes more distributed, accuracy declines. This observation aligns with our theoretical findings. A notable trend is that performance degradation is more marked at higher $\gamma$ values. This suggests a disproportionately negative impact on accuracy when the data is more fragmented, especially if the data quality is poorer.
	
	The cost of adaptation also presents intriguing behavior. We note that the performance gap between adaptive and non-adaptive methods widens with increasing $m$, hinting that distributed data may amplify the challenges of adaptation. When it comes to kernel choice in the context of adaptation, distinct behaviors emerge depending on the specific scenario. For instance, in regimes of high privacy with low data distribution (e.g., $m=1$), the triangular kernel outperforms its Gaussian counterpart in adaptation. Conversely, for larger $m$ values, the Gaussian kernel seems more adept. In scenarios of lower privacy constraints, the triangular kernel consistently outperforms the Gaussian.
	
	From a practical standpoint, the general recommendation is to opt for the Gaussian kernel only when dealing with extensively distributed data (large $m$) coupled with strict privacy requirements. In all other scenarios, the triangular kernel proves to be more effective.
	
	%%%%%%%%%%%%%%%%%%%%%%%%%%%%%%%%%%%%%%%%%%
	\subsection{Real Data}\label{sec:real-data}
	%%%%%%%%%%%%%%%%%%%%%%%%%%%%%%%%%%%%%%%%%%
	
	We now apply our private, distributed, nonparametric classifier to a real dataset to demonstrate its practical merits. We have selected the heart disease dataset from \cite{detrano1989international}, publicly available on the UCI Machine Learning Repository. The primary task is to predict the prevalence of heart disease based on 13 covariates, including demographics (age, sex) and clinical measurements (blood pressure, resting heart rate, cholesterol, chest pain prevalence, etc.). This dataset comprises patient data from four hospitals in Cleveland, Hungary, Switzerland, and Long Beach. The distributed nature of the data, along with the presence of sensitive patient information, makes it an ideal example for applying our distributed differentially private classification algorithm.

	\begin{figure}[h]
		\centering
		\includegraphics[scale = 0.5]{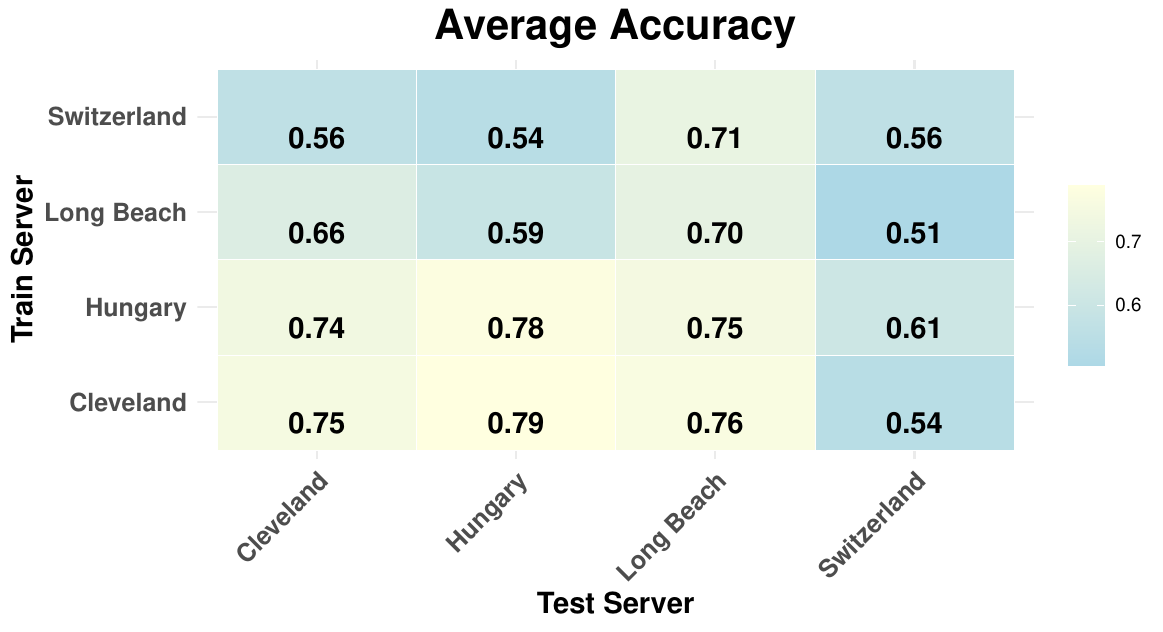}
		\caption{Accuracy for different train-test server pairs}
		\label{fig:server_wise_accuracy}
	\end{figure}

	In order to illustrate the heterogeneity across different servers, we evaluated the performance of our kernel-based classification algorithm in a non-private and non-transfer learning setting. Specifically, we selected 100 random samples from a designated training server and assessed the model's performance not only on a separate test set from the same server but also on data from three additional servers. Figure \ref{fig:server_wise_accuracy} visualizes these results, clearly showing that the same model exhibits varying levels of performance across different servers.
	
	The analysis reveals significant variations in the classification model’s performance, which highlight the inherent heterogeneity in the datasets. Notably, when trained on Cleveland data, the model achieves its highest accuracy on Hungarian data, suggesting some level of similarity between these datasets. In contrast, models trained on Swiss data exhibit markedly lower performance across all servers, including on local data, which points to potential challenges with the Swiss dataset's complexity or representativeness. Meanwhile, the data from Long Beach demonstrates consistent performance across diverse test servers, suggesting its features may possess a higher degree of generalizability. This variability underscores the need for adaptive modeling strategies in distributed systems to manage data heterogeneity and enhance predictive accuracy.
	
	To address these challenges and potentially enhance model performance, we apply our transfer learning models to borrow strength from various sources. We will use the Hungarian hospital as the transfer target, and the rest would be used at source datasets.
	
	%%%%%%%%%%%%%%%%%%%%%%%%%%%%%%%%%%%%%%%%%%
	\subsubsection{Implementation Details}
	%%%%%%%%%%%%%%%%%%%%%%%%%%%%%%%%%%%%%%%%%%
	
	The choice of variables is critical for an illustration of our method, since it is nonparametric and hence extremely susceptible to the curse of dimensionality. Additionally, the data is plagued with missing observations for many of its columns. Among the 13 covariates available, we first picked the 9 predictors which had at least 90\% non-missing observations across all the data. Among the rest, we removed \texttt{fbs} because it had 60\% missingness for Switzerland data. We also removed the categorical variable \texttt{restecg} since its three levels lead to a higher dimensionality, while existing studies suggest its lower importance in predicting heart disease. Our chosen set of variables is therefore of size 7, and consists of \texttt{age}, \texttt{sex}, \texttt{cp} (chest pain type), \texttt{exang} (exercise induced angina, yes or no), \texttt{thalach} (maximum heart rate), \texttt{oldpeak} (ST depression induced by exercise), and \texttt{trestbps} (resting blood pressure). We further excluded the patients who had missing observations for these variables and set aside a test set of size $150$ from our target. This resulted in $n_0=142$ patients for the target (Hungary), while the sources had $n_1=303$, $n_2=141$, and $n_3=116$ patients for the hospitals in Cleveland, Long Beach and Switzerland respectively. 
	
	In order to better place the problem in our theoretical setting, we scale each of the chosen covariates to range between 0 and 0.5. Another important distinction is the centering for our estimators. An initial exploratory analysis for our data reveals that the sources have different degrees of disease prevalence: 47\%, 36\%, 93\%, and 79\% for Cleveland, Hungary, Switzerland, and Long Beach respectively. To alleviate this issue, we re-weight the summands in our kernel estimator by source-specific mean of the disease status.
	
\begin{figure}[h]
		\centering
		\includegraphics[scale = 0.7]{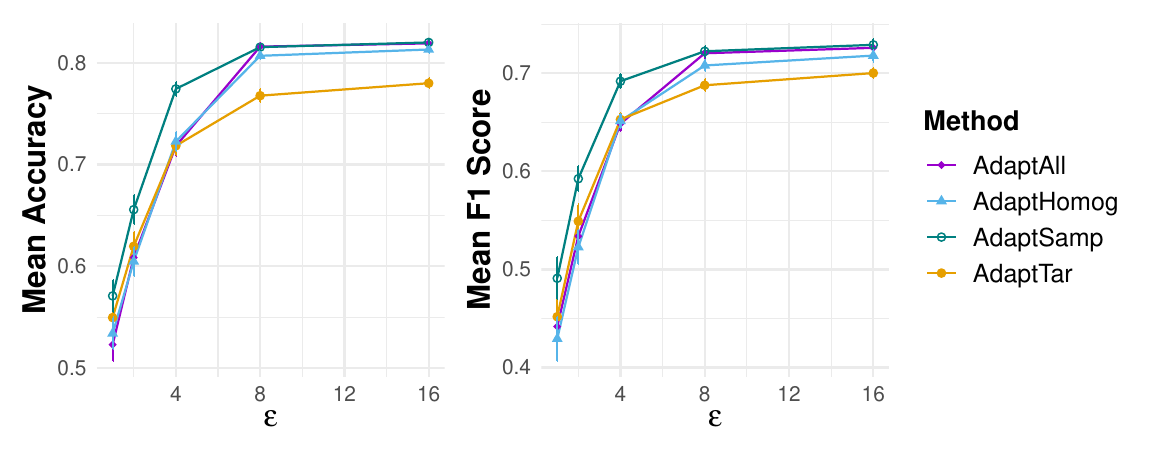}
		\caption{Accuracy and F1 score values for different privacy levels}
		\label{fig:accuracy_and_auc}
	\end{figure}

	As the true transfer parameters are unknown, we opt for the data-adaptive procedures specified in Section~\ref{sec:adaptation}. We will compare the following adaptation variations, each incorporating different amounts of auxiliary information:

	\begin{enumerate}
		\item \textbf{AdaptAll}: We choose the best bandwidth and weights following Section~\ref{sec:adapt-gen}. 
		\item \textbf{AdaptTar}: The best bandwidth is chosen based solely on the target estimator. That is, the weights are zero for each source server.
		%\item \textbf{AdaptEq}: We constrain the weights to be equal for each server and find the best bandwidth.
		\item \textbf{AdaptSamp}: The weights are taken to be proportional to the sample sizes for each server.
		\item \textbf{AdaptHomog}: The adaptation procedure under source homogeneity (see Section~\ref{sec:adapt-homog}) is used to choose the best bandwidth and weights.
	\end{enumerate}
	We compare accuracy of the five different adaptation techniques against the common privacy parameter $\eps$. We also compare the different adaptation procedures on another popular performance metric called the F1 score which is the harmonic mean of precision and recall. Figure~\ref{fig:accuracy_and_auc} shows the average accuracy and F1 score obtained by replicating our procedure $200$ times on different test and train folds. 
	
	It is evident from Figure~\ref{fig:accuracy_and_auc} that both the accuracy and the F1 score values improve for all the classifiers as the privacy requirements become less stringent. Among the five classifiers, \textbf{AdaptTar}, the one depending only on the target shows the poorest performance, thus highlighting the benefits of transfer learning in this problem. On the other hand, the sample size weighted adaptive classifier, denoted \textbf{AdaptSamp} is clearly the best among the five. The all-adaptive classifier, \textbf{AdaptAll}, automatically adjusts to the unknown weights and attains the same accuracy and F1 score values as \textbf{AdaptSamp} for higher values of $\eps$. Finally the \textbf{AdaptHomog} is strictly worse than our estimator possibly indicating the inherent heterogeneity among the source servers. Overall, we find a prediction accuracy of 81\% and an F1 score of 73\% in predicting heart disease from the hospitals dataset, showing the effectiveness of our data-adaptive distributed private transfer learning classification mechanism.
	
	%%%%%%%%%%%%%%%%%%%%%%%%%%%%%%%%%%%%%%%%%%
	\section{Discussion}\label{sec:discussion}
	%%%%%%%%%%%%%%%%%%%%%%%%%%%%%%%%%%%%%%%%%%
	
	In this paper, we establish the minimax misclassification rate in a heterogeneous distributed setting with varying sample sizes, privacy parameters, and data distributions across servers under the posterior drift model. Our results precisely characterize the effects of privacy constraints, source sample sizes, and target sample size.
	
	We rigorously quantify the trade-offs between data heterogeneity and variations in privacy budgets. By analyzing these trade-offs through minimax optimality, we clarify how differences in data distribution across servers impact the overall learning process and model performance. Our findings highlight the critical balance needed between maintaining privacy and effectively leveraging distributed data. The construction of the minimax and data-driven adaptive classifiers addresses excess risk and encapsulates the inherent trade-offs introduced by differential privacy. The impact of privacy constraints can be mitigated, to some extent, by leveraging larger datasets. The nuanced behavior of these classifiers, specifically how they scale with changes in the privacy constraints and sample sizes, offers a promising avenue for optimizing distributed learning systems. Our results are robust and theoretically justified within the defined heterogeneous settings and the posterior drift model. 
	
	%Following our results, a number of intriguing related problems can be considered. From a technical standpoint, it is of interest to extend this result to the case where the regression functions $\eta_j(\cdot)$ are $\beta$-H{\"o}lder with $\beta>1$ (see Remark~\ref{rem:beta<1}). This might be possible using kernels with vanishing moments of order $\floor{\beta}$ along with additional transfer exponent assumptions on the higher order derivatives of $\eta_j$. In addition to the complications of transfer learning when $\beta>1$, it is challenging to find a suitable function privatization method akin to the one used in this paper. Going beyond kernel based estimators, \cite{cai2021transfer} have suggested using weighted nearest neighbor estimator for transfer learning when $\beta>2$ --- finding an optimal privatization method for such nearest neighbor based methods would require fundamentally different techniques from the current work, and is a potential direction for future research.
	
	Motivated by several practical requirements, it would be of significant interest to consider the transfer learning problem in other models with distributed differential privacy constraints. One natural direction is to consider distributed classification under the covariate shift model with differential privacy constraints. Such an analysis could provide additional insights and methods for a range of potential applications. Another future direction is to study distributed classification under parametric models such as logistic regression. We leave these intriguing problems for future research.

\bigskip  
\noindent \textbf{Funding}: The research was supported in part by NIH grants R01-GM123056 and R01-GM129781.

	\bibliographystyle{plainnat}
	\bibliography{references.bib}

\end{document}